\documentclass[10pt]{amsart}
\usepackage{amsfonts}
\usepackage{amsmath}
\usepackage{amsthm}

\newtheorem{lem}{Lemma}

\newtheorem{prop}{Proposition}

\theoremstyle{definition}
\newtheorem{defn}{Definition}
\theoremstyle{definition}
\newtheorem{thm}{Theorem}
\newtheorem*{conj}{Conjecture}

\newenvironment{pf}{\proof}{\endproof}

\theoremstyle{remark}

\numberwithin{equation}{section}

\begin{document}

\newcommand{\thmref}[1]{Theorem~\ref{#1}}
\newcommand{\secref}[1]{Sect.~\ref{#1}}
\newcommand{\lemref}[1]{Lemma~\ref{#1}}
\newcommand{\propref}[1]{Proposition~\ref{#1}}
\newcommand{\corref}[1]{Corollary~\ref{#1}}
\newcommand{\remref}[1]{Remark~\ref{#1}}
\newcommand{\er}[1]{(\ref{#1})}
\newcommand{\nc}{\newcommand}
\nc{\goth}{\mathfrak}
\renewcommand{\Bbb}{\mathbb}
\nc{\Xp}[1]{X^+(#1)}
\nc{\Xm}[1]{X^-(#1)}
\nc{\on}{\operatorname}
\nc{\ch}{\mbox{ch}}
\nc{\Z}{{\Bbb Z}}
\nc{\J}{{\mathcal J}}
\nc{\C}{{\Bbb C}}
\nc{\Q}{{\Bbb Q}}
\renewcommand{\P}{{\mathcal P}}
\nc{\N}{{\Bbb N}}
\nc{\ov}{\overline}
\nc{\pone}{{\Bbb C}{\Bbb P}^1}
\nc{\pa}{\partial}
\def\H{\mathcal H}
\def\L{\mathcal L}
\nc{\F}{{\mathcal F}}
\nc{\Sym}{{\goth S}}
\nc{\A}{{\mathcal A}}
\nc{\arr}{\rightarrow}
\nc{\larr}{\longrightarrow}
\nc{\al}{\alpha}
\nc{\ri}{\rangle}
\nc{\lef}{\langle}
\nc{\W}{{\mathcal W}}
\nc{\la}{\lambda}
\nc{\ep}{\epsilon}
\nc{\su}{\widehat{\goth{sl}}_2}
\nc{\asltwo}{\widehat{\goth{sl}}_2}
\nc{\sw}{\goth{sl}}
\nc{\g}{\goth{g}}
\nc{\h}{\goth{h}}
\nc{\n}{\got{n}}
\nc{\ab}{\goth{a}}
\nc{\G}{\widehat{\g}}
\nc{\De}{\Delta}
\nc{\gt}{\widetilde{\g}}
\nc{\Ga}{\Gamma}
\nc{\one}{{\bold 1}}
\nc{\hh}{\widehat{\h}}
\nc{\z}{{\goth Z}}
\nc{\zz}{{\mathcal Z}}
\nc{\Hh}{{\mathcal H}}
\nc{\qp}{q^{\frac{k}{2}}}
\nc{\qm}{q^{-\frac{k}{2}}}
\nc{\La}{\Lambda}
\nc{\wt}{\widetilde}
\nc{\qn}{\frac{[m]_q^2}{[2m]_q}}
\nc{\cri}{_{\on{cr}}}
\nc{\sun}{\widehat{\sw}_N}
\nc{\HH}{{\mathcal H}_q(\sw_N)}
\nc{\ca}{\wt{{\mathcal A}}_{h,k}(\sw_2)}
\nc{\si}{\sigma}
\nc{\gl}{\goth{g}\goth{l}_N}
\nc{\el}{\ell}
\nc{\s}{t}
\nc{\tN}{\theta_{p^N}}
\nc{\ds}{\displaystyle}
\nc{\Dp}{D_{p^{-1}}}
\nc{\Dq}{D_{q^{-1}}}
\nc{\tL}{{\bold L}}
\nc{\tP}{{\bold P}}
\nc{\tA}{{\bold A}}
\nc{\tU}{{\bold U}}
\nc{\beq}{\begin{equation}}
\nc{\eeq}{\end{equation}}
\nc{\PN}{{\mathcal P}_{q,N}}
\nc{\RN}{{\mathcal R}_{q,N}}
\nc{\LN}{{\mathcal L}_{q,N}}
\nc{\Uq}{{\mathcal U}_{q}}
\nc{\bi}{\bibitem}
\nc{\uqatwo}{{U_q}(\su)}
\nc{\uaqtwo}{{U_q}(\su)}
\nc{\uaqtwoatone}{{U_{q,1}}(\su)}
\nc{\uqone}{{U_{q,1}}(\su)}
\nc{\uqatwoatone}{{U_{q,1}}(\su)}
\nc{\uqtwo}{U_q(\goth{sl}_2)}
\nc{\dij}{\delta_{ij}}
\nc{\divei}{E_i^{(n)}}
\nc{\divfi}{F_i^{(n)}}
\nc{\Lzero}{\Lambda_0}
\nc{\Lone}{\Lambda_1}
\nc{\ve}{\varepsilon}
\nc{\phioneminusi}{\Phi^{(1-i,i)}}
\nc{\phioneminusistar}{\Phi^{* (1-i,i)}}
\nc{\phii}{\Phi^{(i,1-i)}}
\nc{\Li}{\Lambda_i}
\nc{\Loneminusi}{\Lambda_{1-i}}
\nc{\vtimesz}{v_\ve \otimes z^m}
\nc{\xkpm}{x_k^{\pm}}
\nc{\xlpm}{x_l^{\pm}}
\nc{\xnpm}{x_n^{\pm}}
\nc{\xkp}{x_k^{+}}
\nc{\xkm}{x_k^{-}}
\nc{\xlp}{x_l^{+}}
\nc{\xlm}{x_l^{-}}
\nc{\xzp}{x_0^{+}}
\nc{\xzm}{x_0^{-}}
\nc{\xom}{x_1^{-}}
\nc{\xmp}{x_{-1}^{+}}
\nc{\ppm}{\Phi_{\pm}}
\nc{\pP}{\Phi_{+}}
\nc{\pM}{\Phi_{-}}
\nc{\bppm}{\ov{\Phi}_{\pm}}
\nc{\bp}{\ov{\Phi}}
\nc{\bpP}{\ov{\Phi}_{+}}
\nc{\bpM}{\ov{\Phi}_{-}}
\nc{\vi}{|\Lambda_i\rangle}
\nc{\voneminusi}{|\Lambda_{1-i}\rangle}
\nc{\vone}{|\Lambda_1\rangle}
\nc{\vzero}{|\Lambda_0\rangle}

\nc{\expminusaminus}[1]{\exp\bigl(- \sum_{n = 1}^{\infty}
  \frac{q^{{#1}n/2}}{n} a_{-n} z^{n}\bigr)}
\nc{\expminusaplus}[1]{\exp\bigl(- \sum_{n = 1}^{\infty}
  \frac{q^{{#1}n/2}}{n} a_{n} z^{-n}\bigr)}
\nc{\expplusaminus}[1]{\exp\bigl( \sum_{n = 1}^{\infty}
  \frac{q^{{#1}n/2}}{n} a_{-n} z^{n}\bigr)}
\nc{\expplusaplus}[1]{\exp\bigl( \sum_{n = 1}^{\infty}
  \frac{q^{{#1}n/2}}{n} a_{n} z^{-n}\bigr)}

\nc{\expminusaminuszsame}[1]{\exp\bigl(- \sum_{n = 1}^{\infty}
  \frac{q^{{#1}n/2}}{n} a_{-n} z^{-n}\bigr)}
\nc{\expminusapluszsame}[1]{\exp\bigl(- \sum_{n = 1}^{\infty}
  \frac{q^{{#1}n/2}}{n} a_{n} z^{n}\bigr)}
\nc{\expplusaminuszsame}[1]{\exp\bigl( \sum_{n = 1}^{\infty}
  \frac{q^{{#1}n/2}}{n} a_{-n} z^{-n}\bigr)}
\nc{\expplusapluszsame}[1]{\exp\bigl( \sum_{n = 1}^{\infty}
  \frac{q^{{#1}n/2}}{n} a_{n} z^{n}\bigr)}

\nc{\macsum}[4]{\sum_{l(\lambda) \le N} P_\lambda(#1;
  q^{#3},q^{#4})\P_{\lambda'}({#2};q^{#4},q^{#3})}
\nc{\macsumq}[4]{\sum_{l(\lambda) \le N} P_\lambda(#1;
  q^{#3},q^{#4})Q_{\lambda}({#2};q^{#3},q^{#4})}
\nc{\macsummu}[4]{\sum_{l(\mu) \le N} P_\mu(#1;
  q^{#3},q^{#4})\P_{\mu'}({#2};q^{#4},q^{#3})}
\nc{\macsumqmu}[4]{\sum_{l(\mu) \le N} P_\mu(#1;
  q^{#3},q^{#4})Q_{\mu}({#2};q^{#3},q^{#4})}

\nc{\Nplusmacsum}[4]{\sum_{l(\lambda) \le N+1} P_\lambda(#1;
  q^{#3},q^{#4})\P_{\lambda'}({#2};q^{#4},q^{#3})}
\nc{\Nplusmacsumq}[4]{\sum_{l(\lambda) \le N+1} P_\lambda(#1;
  q^{#3},q^{#4})Q_{\lambda}({#2};q^{#3},q^{#4})}
\nc{\Nplusmacsummu}[4]{\sum_{l(\mu) \le N+1} P_\mu(#1;
  q^{#3},q^{#4})\P_{\mu'}({#2};q^{#4},q^{#3})}
\nc{\Nplusmacsumqmu}[4]{\sum_{l(\mu) \le N+1} P_\mu(#1;
  q^{#3},q^{#4})Q_{\mu}({#2};q^{#3},q^{#4})}

\nc{\macsump}[4]{\sum_{l(\lambda) \le N} P_\lambda(#1;
  q^{#3},q^{#4})P_{\lambda}({#2};q^{#4},q^{#3})}

\nc{\eali}[1]{e^{#1 \alpha}\vi}
\nc{\eal}[1]{e^{#1 \alpha}}
\nc{\ealoneminusi}[1]{e^{#1 \alpha}\voneminusi}
\nc{\expplusaminusz}[1]{\exp\bigl( \sum_{n = 1}^{\infty}
  \frac{q^{{#1}n/2}}{n} a_{-n} \sum_{j=1}^N z_j^{n}\bigr)}

\nc{\expminusaminusz}[1]{\exp\bigl(- \sum_{n = 1}^{\infty}
  \frac{q^{{#1}n/2}}{n} a_{-n} \sum_{j=1}^N z_j^{n}\bigr)}

\nc{\prodz}{\prod_{j = 1}^N z_j^{\lfloor \frac{N-j+1+i}{2}\rfloor }}
\nc{\prodzminus}{\prod_{j = 1}^N z_j^{-\lfloor \frac{N-j+1-i}{2}\rfloor }}
\nc{\prodxi}{\prod_{i < j} \xi(\frac{z_j}{z_i})}
\nc{\minusonepower}{(-1)^{i + \lfloor (N -i)/2\rfloor \lfloor (N-1-i)/2\rfloor
}}

\nc{\xkplusonepm}{x_{k+1}^{\pm}}
\nc{\xlplusonepm}{x_{l+1}^{\pm}}
\nc{\xkpluslpm}{x_{k+l}^{\pm}}
\nc{\hxzp}{\hat x_0^+}
\nc{\hxzm}{\hat x_0^-}
\nc{\hxmp}{\hat x_{-1}^+}
\nc{\hxpm}{\hat x_{1}^-}
\nc{\fracqqi}[1]{\frac{#1}{q+q^{-1}}}
\nc{\ahat}{(-1)^{n-1} a}
\nc{\ahatn}{(-1)^{n} a}
\nc{\td}{\vdots}
\nc{\symfun}{\Lambda^{\goth S}}
\title{Canonical Basis and Macdonald Polynomials}
\dedicatory{To Howard Garland on his 60th birthday}
\author{Jonathan Beck} \
\address{Jonathan Beck, University of Chicago}
\author{Igor B. Frenkel}
\address{Igor B. Frenkel, Yale University}
\author{Naihuan Jing}
\address{Naihuan Jing, North Carolina State University}
\begin{abstract}
  In the basic representation of $\uqatwo$ realized via the algebra of
  symmetric functions we compare the canonical basis with the basis of
  Macdonald polynomials with $t=q^2$.  We show that the Macdonald
  polynomials are invariant with respect to the bar involution defined
  abstractly on the representations of quantum groups.  We also prove
  that the Macdonald scalar product coincides with the abstract
  Kashiwara form.   This implies, in particular, that the Macdonald
  polynomials form an intermediate basis between the canonical basis and
  the dual canonical basis, and the coefficients of the transition
  matrix are necessarily bar invariant.   We also discuss the positivity
  and integrality of these coefficients.  For level $k$, we expect a
  similar relation between the canonical basis and Macdonald polynomials
  with $q^2 = t^{k}.$
\end{abstract}
\maketitle
\pagestyle{plain}


\date{Advances in Math., accepted}

\maketitle

\centerline{Advances in Math., to appear}
\section{Introduction.}

Since G. Lusztig \cite{Lu} and M. Kashiwara \cite{Ka} introduced
canonical bases of highest weight representations of Kac--Moody
algebras, there has been a significant effort to find their explicit
algebraic description.  The problem becomes especially intriguing for
the simplest affine Lie algebra $\su.$ It is known that the two basic
representations of $\su$ have a realization in the space of symmetric
functions in infinitely many variables $\C[a_{-1}, a_{-2}, \dots ],$
tensored with the extra space (of ``sectors'') $\oplus \C e_n, n \in \Z$
or $n \in \Z + \frac{1}2$ \cite{FK}, \cite{S}, and this construction
admits a $q$ deformation \cite{FJ}.  Therefore the canonical basis
yields a class of symmetric polynomials depending on a parameter $q$.
Moreover, integrable representations of level $k$ for $\su$, having in
general more complicated structure, contain natural subspaces generated
by $e(z)^k, f(z)^k,$ that also have a symmetric space realization
\cite{LP}.  Again, this construction admits a $q$--deformation \cite{J3,
  DF} and leads to a two parameter family of symmetric functions.

A few years before the discovery of the canonical basis, I. Macdonald
\cite{Ma} found a remarkable two parameter family of symmetric functions
that includes as a special case practically all known classical
symmetric functions.  A representation theoretic intepretation of
Macdonald polynomials in terms of certain spherical functions for the
quantum groups $U_q({\frak sl}_n)$ was given in \cite{EK}. A vertex
operator approach to Hall-Littlewood and some Macdonald polynomials 
can be found in \cite{J1, J2}.
Macdonald
polynomials can also be viewed as a basis of the symmetric space
$\C[a_{-1}, a_{-2}, \dots]$, and it is natural to question
its relation to the canonical basis for the quantum affine algebra
$\uqatwo.$

To compare the two bases we recall a characterization of the canonical
basis given by Kashiwara \cite{Ka}.  The elements of the canonical basis
are determined up to a sign by the following properties: 1) they are
invariant with respect to a bar involution, an operation defined on
any highest weight representation of the quantum affine algebra
associated to a Kac--Moody algebra, 2) they are orthogonal modulo
$q^{-1}\Z[q^{-1}]$ with respect to a bilinear form, defined by Kashiwara
on any highest weight representation and, 3) they belong to the lattice
of divided powers of the quantum group over $\Z[q,q^{-1}].$ It turns out
that the basis of Macdonald polynomials satisfy similar (and in a
certain sense simpler) properties than the ones that characterize the
canonical basis.  We show that the Macdonald polynomials form an
``intermediate basis'' between the dual canonical basis and the
canonical basis with respect to the Kashiwara form.  We also show that
the transition matrix between the basis of Macdonald polynomials and the
canonical basis (after a minor rescaling) is bar invariant, integral,
and independent of the lattice point $e_n.$ We conjecture its
positivity.

In this paper we only consider the level $1$ representations, and this
corresponds to the relation $t = q^2$ (or $q = t^2$) between the two
parameters of the Macdonald polynomials.  An arbitrary level $k$
representation will be considered in a sequel to this paper.  In Section
$2$, we recall basic facts about the quantum affine algebra $\uqatwo$ and
its representations.  We also give a definition and characterization of
the canonical basis.  In Section $3$, we give the loop--like presentation
of the quantum affine algebra and realization of basic representations
in the space of symmetric functions.  We then introduce certain
intertwining operators for $\uqatwo$ and express them via quantum vertex
operators.  In Section $4$, we recall the definition and properties of
Macdonald polynomials and identify their generating functions as
``one--half'' of quantum vertex operators.

In the next three sections we study the properties of the Macdonald
polynomials that correspond to the three characteristic properties of
the canonical basis, and establish how the two structures are related.
It is the quantum vertex operator, viewed as an intertwining operator
for $\uqatwo$, and as a generating function for Macdonald polynomials,
which provides a bridge between the two theories.  In Section 5, we
prove the invariance of Macdonald polynomials (rescaled by appropriate
powers of $q$, depending on sector) under bar conjugation.  In Section
6, we establish the coincidence of the Macdonald and Kashiwara forms.
This implies that the Macdonald polynomials are orthogonal with respect
to the Kashiwara form, and that the transition matrix between the dual
canonical basis and the canonical basis for a fixed weight admits a
decomposition
$$ A(q)^t D(q) A(q),$$ where $A(q) = \{a_{\lambda,\mu}\}$ is a bar
invariant matrix, and $D(q)$ is diagonal and bar invariant.  The matrix
$A(q)$ is precisely the transition matrix between the canonical basis
and the basis of Macdonald polynomials.  The coincidence of the
Macdonald and Kashiwara forms implies that the matrix $A(q)$ depends
trivially on the choice of a sector.  In Section 7, we show that the
$\Z[q,q^{-1}]$--lattice spanned by the dual integral Macdonald
polynomials (rescaled Macdonald polynomials introduced in \cite{Ma})
form a lattice in the basic representation which is invariant under the
divided powers action.  Finally, in Section 8, we show that the integral
Macdonald polynomials belong to the lattice of divided powers of the
quantum group over $\Z[q,q^{-1}].$ This implies the integrality of the
rescaled matrix coefficients, and we also conjecture their positivity.

Our results on the relation between canonical basis and Macdonald
polynomials open new avenues of research in both subjects.  On one hand,
Macdonald polynomials suggest a definition of ``intermediate canonical
basis'' in any integrable highest weight representation of Kac--Moody
algebras, with the property of bar--invariance and orthogonality with
respect to Kashiwara's form.  On the other hand, the canonical basis and
the dual canonical basis in highest weight representations of $\uqatwo$
yield a pair of dual bases of symmetric functions with respect to
Macdonald's form.  These bases, depending on two parameters, have a
natural representation theoretic origin and should have significance in
the theory of symmetric functions.  Further understanding of the
relation between the two theories, most importantly the coefficients of
the transition matrix $A(q)$, should reveal new layers of structure
behind these profound subjects.


\section{Preliminary material.}
\subsection{Quantum algebras}
Let $\A = \Z[q, q^{-1}].$   For $n \in \N$ we define
\begin{equation*} [n] = \frac{q^n
    - q^{-n}}{q - q^{-1}}, \ [n]! = [n] [n-1] \dots [1], \ [0]! = 1, \
\end{equation*}
Let the (derived) {\em quantum affine algebra} $\uqatwo$ be generated over $\C(q)$ by
$E_i, \ F_i, K_i^{\pm 1},\ i = 0,1,$ with relations:
\begin{align*}
  & K_iE_jK_i^{-1} = q^{a_{ij}} E_j,\ \ K_iF_jK_i^{-1} = q^{-a_{ij}}
 F_j, \ \ [E_i,F_j]
  = \dij \frac{K_i - K_i^{-1}}{q - q^{-1}}, \\ & E_i^{(3)}E_j -
  [3]E_i^{(2)}E_jE_i + [3]E_iE_jE_i^{(2)} + E_jE_i^{(3)}=0 ,\ i \neq j, \\
  & F_i^{(3)}F_j - [3]F_i^{(2)}F_jF_i + [3]F_iF_jF_i^{(2)} +
  F_jF_i^{(3)}=0, \ i \neq j,
\end{align*}
where $a_{ii} = -a_{ij} = 2$ for $i \neq j.$  

Let $E_i^{(n)} = \frac{E_i^n}{[n]!}, \ F_i^{(n)} = \frac{F_i^n}{[n]!}.$
The quantum algebras have an $\A$--form $U_\A,$ generated over $\A$ by
$\divei, \ \divfi,$ and $K_i^{\pm 1}$.
We define a Hopf algebra structure  on $\uqatwo$ with coproduct $\De$,
antipode $S$, and counit $\eta$ given  by:
\begin{align*} \label{coproduct}
&\De(E_i) = K_i \otimes E_i + E_i \otimes 1, \
\De(F_i) = 1 \otimes F_i + F_i \otimes K_i^{-1}, \
\De(K_i) = K_i \otimes K_i.  \\
& S(E_i) = - K_i^{-1} E_i, \ S(F_i) = -F_i K_i, \ S(K_i) = K_i^{-1}. \\
&\eta(E_i) = 0, \ \eta(F_i) = 0, \ \eta(K_i) = 1.
\end{align*}
This coproduct structure is not unique.  Given an (anti)--automorphism
of the algebra $\sigma$, the map $(\sigma \otimes \sigma) \circ \Delta
\circ \sigma^{-1}$ is also a coproduct.
\subsection{Representations}
 The free abelian group $X$ on
$\{ \Lzero, \Lone, \delta \}$
is called the {\em weight
lattice} with {\em fundamental weights} $\Lambda_i, \ i = 0, 1.$
Define the {\em simple roots} $\alpha_0, \alpha_1$ by
$$\alpha_0 + \alpha_1 = \delta, \ \Lambda_1 = \Lambda_0 +
(\alpha_1/2).$$

For $\lambda \in X^+,$ let $V = V(\lambda)$ be the irreducible
integrable representation of $\uqatwo$ with highest weight $\lambda.$ We
denote by $V(\lambda)_\A$ the image of the canonical map $U_\A
\rightarrow V_\lambda.$ $V$ is graded by $X^+$ and decomposes as $V =
\oplus_{\nu} V_\nu.$ An element $x \in V_\nu$ is said to be of {\em
  weight $\nu$}.  Denote by $\hat V$ the completion of $V$ with respect
to this homogeneous grading.

In addition to the above highest weight representations, we will also
make use of the two dimensional evaluation representation \cite{Jb} of
$\uqatwo$.  Let $\uqtwo$ be the Hopf subalgebra of $\uqatwo$ generated
over $\Q(q)$ by $E_1$, $F_1$ and $K_1^{\pm1}$.  Let $V_1 = {\C}v_+
\oplus {\C}v_-.$ Define the standard two dimensional representation of
$\uqtwo$ by:
\begin{equation}
E_1 v_- = v_+, \  E_1 v_+ = 0, \ F_1 v_+ = v_-, \ F_1 v_- = 0, K_1 v_{\pm} =
q^{\pm 1} v_\pm.
\end{equation}
Let $z$ be an indeterminate and consider $V_z = V \otimes \C[z,z^{-1}].$
We define a $\uqatwo$ action as follows:
\begin{align*}
& E_0(v_\ve \otimes z^m) = (F_1 v_\ve) \otimes z^{m+1}, \ E_1(\vtimesz)
= (E_1 v_{\ve}) \otimes z^{m}, \\
& F_0(v_\ve \otimes z^m) = (E_1 v_\ve) \otimes z^{m-1}, \ F_1(\vtimesz)
= (F_1 v_{\ve}) \otimes z^{m}, \\
& K_1(\vtimesz) = (K_1 v_\ve) \otimes z^m, \
 K_0 (\vtimesz) = K_1^{-1}(\vtimesz),
\end{align*}
where $\ve = \pm.$

\begin{defn} The algebra {\em $\uaqtwo$ at level $k$}
is defined as
\begin{equation*} U_{q,k}({\asltwo}) = \frac{\uqatwo}{(C-q^k) \uqatwo}.
\end{equation*}
\end{defn}
There
is a natural surjection from any representation $V$ of $\uqatwo$
to a corresponding representation of the algebra $U_{q,k}({\asltwo}).$

\subsection{Canonical Basis}
Introduce the $\Q(q)$--linear anti--involution $\rho$, and the
$\Q$--linear involution $^{\overline{\ \ }}$
of $\uaqtwo$ by:
\begin{align}
\label{rho}
&\rho(E_i) = q K_i F_i, \ \rho(F_i) = q K_i^{-1} E_i, \ \rho(K_i) = K_i, \\
& \label{bar} \overline{E_i} = E_i, \ \overline{F_i} = F_i, \ \overline{K_i} =
K_i^{-1}, \overline q = q^{-1}.
\end{align}
\begin{prop} \cite{Ka}
For $\lambda \in X^+$, let $V = V(\lambda).$
There is a unique bilinear form $(\ ,\ ) : V \times V \rightarrow \Q(q)$
such that:
\begin{align}
\label{highestweight}
& (v_\lambda, v_\lambda) = 1, \\
\label{rhochange} & (ux, y) =
(x, \rho(u)y)  \text{ for all } x, \ y \in V \text{ and }
  u \in \Uq.
\end{align}
This bilinear form is symmetric.  If  $x \in V_\nu, y \in
V_{\nu'}$  with  $\nu \neq \nu'$,  then  $(x,y) = 0.$
\end{prop}

We cite the following characterization of the canonical basis given by
Kashiwara:
\begin{thm} \cite{Ka} \label{Lusztigtheorem}
 Let $b \in V(\lambda)_{\A}$
Then either $b$ or $-b$ is in the canonical basis if and only if
\begin{enumerate}
\item $\bar b = b$, and
\item $(b,b') = \delta_{b,b'}$ mod $q^{-1} \Z[q^{-1}].$
\end{enumerate}
\end{thm}


\section{Realization of basic representation and intertwiners}

\subsection{Loop--like realization}
The algebra $\uaqtwo$ has another set of loop--like generators.
\begin{prop} \cite{Dr, Be} ,  $\uqatwo$ is isomorphic to the algebra
generated over
$\Q(q)$ on the generators $a_n\ (n \in \Z \setminus \{0\}), \ \xkpm\ (k \in
\Z),\  C^{\pm \frac{1}2},\ K^{\pm 1}$ and the following relations:
\begin{align}
& \label{aas}
[a_n, a_m] = \delta_{m, -n} \frac{n}{[2n]} \frac{C^n -  C^{-n}}{q -
  q^{-1}}  \\
& K a_n K^{-1} = a_n, \ K \xkpm K^{-1} = q^{\pm 2} \xkpm, \notag \\
& [a_l, \xkpm] = \pm  C^{\mp |l|/2} \xkpluslpm, \label{aandx} \\
& \xkplusonepm \xlpm  - q^{\pm 2} \xlpm \xkplusonepm = q^{\pm 2}  \xkpm
\xlplusonepm - \xlplusonepm \xkpm, \notag \\
& [\xkp, \xlm] = \frac{C^{(k-l)/2} \psi_{k+l} - C^{(l-k)/2} \varphi_{k+l}}{q -
q^{-1}}. \notag
\end{align}
\end{prop}
where $C$ is a central element and
\begin{align}
& \sum_{k=0}^\infty \psi_k z^{-k} = K
\exp\bigl((q-q^{-1})\sum_{k=1}^\infty \frac{[2k]}{k} a_k z^{-k}\bigr), \\
& \sum_{k=0}^\infty \varphi_{-k} z^{k} = K^{-1}
\exp\bigl(-(q-q^{-1})\sum_{k=1}^\infty  \frac{[2k]}{k} a_{-k} z^{k}\bigr).
\end{align}
The isomorphism is determined by mapping the respective generators as
follows:
\begin{align*}
&K_1 \mapsto  K, \ E_1 \mapsto x_0^+, \ F_1 \mapsto x_0^-, \\
&K_0 \mapsto CK^{-1}, \ E_0K_1 \mapsto x_1^{-}, K_1^{-1} F_0 \mapsto
x_{-1}^+.
\end{align*}
The isomorphism is explained completely in terms of a braid group
action on $\uaqtwo$ in \cite{Be}.
  Forming generating series from the loop--like
generators by
$$X^\pm(z) = \sum_{n \in \Z} \xnpm z^{-n-1},$$ the defining relations are
written as:
\begin{align}
&[a_k, X^\pm(z)] = \pm C^{\mp |k|/2} z^k X^\pm(z), \\
&(z - q^{\pm 2} w) X^\pm(z) X^\pm(w) + (w - q^{\pm 2}z) X^\pm(w)
X^\pm(z) = 0, \\
\begin{split}
&[X^+(z), X^-(w)] = K \exp\big[{(q - q^{-1}) \sum_{k=1}^\infty \frac{[2k]}{k}
  a_k C^{k/2}
  z^{-k}}\big] \frac{\delta(z/Cw)}{(q-q^{-1})zw} \\
& \quad \qquad - K^{-1} \exp\big[{-(q-q^{-1})
\sum_{k =1}^\infty \frac{[2k]}{k} a_{-k} C^{k/2} z^k}\big]
\frac{\delta(Cz/w)}{(q-q^{-1})zw},
\end{split}  \\ \notag
\end{align} where as usual $\delta(z) = \sum_{n \in \Z} z^n.$

\subsection{Basic Representation}
Now we consider
$\uqatwoatone.$   For $i = 0, 1,$ let $V(\Lambda_i)$
be the unique irreducible highest weight representation  with
highest weight $\Lambda_i.$

Let
\begin{equation} \label{basic} V'(\Li) = \C[a_{-n}, n > 0]
  \oplus (\oplus_{n \in \Z} \C e^{\Li + n \alpha}), \ i = 0, 1.
\end{equation} For $a_n (n \neq 0), \ e^\alpha, \ \partial$
define an action on $V'(\Li)$ as follows:
\begin{align} \label{action}
  a_n(f \otimes e^\beta) & = a_n f \otimes e^\beta, \text{ if } n <
  0, \\
                           & = [a_n,f] \otimes e^\beta, \text{ if } n
                           > 0, \notag \\
   e^\alpha (f \otimes e^\beta) & = f \otimes e^{\beta + \alpha},
   \notag \\
   \partial (f \otimes e^\beta) & = (\alpha, \beta) f \otimes
   e^\beta,  \notag
\end{align}
where $f \in \C[a_{-n}, \ n>0]$ and $\beta = \Li + n \alpha.$

\begin{thm} \cite{FJ}
  The representations $V'(\Li)$ and $V(\Li)$ are isomorphic.  The action
of the loop-like generators on $V'(\Li)$ are determined by the
following relations:
\begin{align}
&K = q^\partial , \ C = q, \notag \\
\begin{split}
&X^\pm(z) =   \exp(\pm \sum_{n=1}^\infty \frac{a_{-n} q^{\mp n/2}(q^n  +
    q^{-n})}{n} z^n) \times \\
& \quad \quad \exp(\mp \sum_{n=1}^\infty \frac{a_{n} q^{\mp n/2}(q^n  +
    q^{-n})}{n} z^{-n})e^{\pm \alpha} z^{\pm \partial}.
\end{split} \label{loopaction}
\end{align}
\end{thm}

\subsection{Quantum Vertex Operators}
The tensor product of any two representations of $\uaqtwo$ is defined
via a coproduct.
Vertex operators of type I are
intertwiners (i.e. $\uqatwo$ module homomorphisms) of the following
form:
\begin{equation} \label{phioneminusi}
\phioneminusi(z) : V(\Li) \rightarrow V(\Loneminusi) \otimes V_z.
\end{equation}
The precise meaning of this map is as follows:
\begin{align*}
&\phioneminusi(z) = \sum_{\ve = +, -} \phioneminusi_{\ve}(z) \otimes
v_\ve, \\
& \text{ where }
\phioneminusi_\ve(z) = \sum_{n \in \Z} \phioneminusi_{\ve, n}
z^{-n}.
\end{align*}
Each $\phioneminusi_{\ve, n}$ is defined to be a linear map
\begin{equation*}
\phioneminusi_{\ve, n} : V(\Li) \rightarrow V(\Loneminusi)
\end{equation*}
which intertwines $\uqatwo$ in the sense that:
\begin{equation} \label{interproperty}
\sum_{\ve} \phioneminusi_{\ve, n} x v \otimes (v_\ve \otimes z^{-n})
=  \Delta(x) \{\sum_{\ve} \phioneminusi_{\ve, n} v \otimes (v_\ve
\otimes z^{-n})\}
\end{equation}
for all $x \in \uaqtwo$ and $v \in V(\Li).$
Equivalently
\begin{equation*}
\Delta(x) \circ \Phi(z) = \Phi(z) \circ x \quad \text{for } x \in
\uqatwo.
\end{equation*}
Let $\vzero \in V(\Lambda_0)$ and $\vone \in V(\Lambda_1)$ be
the respective highest weight vectors.
 These operators are further normalized so that
\begin{equation} \label{normalization}
\begin{split} &
\Phi_-(z) \vzero = \vone \otimes v_- + \text{lower weight terms in the
  first
component},  \\
&\Phi_+(z) \vone = \vzero \otimes v_+ + \text{lower weight terms in the
  first
component}.
\end{split}
\end{equation}
Composition of vertex operators is natural. For example,
\begin{equation}
\Phi^{(0,1)}(z_1) \Phi^{(1,0)}(z_2) : V(\Lambda_0) \rightarrow
  V(\Lambda_0) \otimes V_{z_1} \otimes V_{z_2}.
\end{equation}
Defining $\H = V(\Lambda_0) \oplus V(\Lambda_1)$ we have
\begin{equation*}
\Phi := \phii(z) \oplus \phioneminusi(z) : \H \rightarrow \H \otimes
V_z.
\end{equation*}
Define the ``dual'' vertex operators as the unique
intertwiners of the form:
\begin{equation} \label{phioneminusistar}
\phioneminusistar(z) : V(\Loneminusi) \otimes V_z \rightarrow  V(\Li),
\end{equation}
with the components
\begin{equation*}
   \phioneminusistar_\ve(z) v = \phioneminusistar(z) (v
     \otimes v_\ve),
\end{equation*}
and the normalization given by:
\begin{equation*}
\begin{split}
\Phi^*_+(z) \vzero = \vone \otimes v_+ + \text{lower weight terms in the
  first
component},  \\
\Phi^*_-(z) \vone = \vzero \otimes v_- + \text{lower weight terms in the
  first
component}.
\end{split}
\end{equation*}

\begin{thm} \cite{FR} Fix a coproduct $\Delta$ of $\uqatwo$.
The vertex operator $\phioneminusi(z)$ exists
and is uniquely determined by highest weight normalization.
The product of vertex operators
\begin{equation} \Phi_{\ve_1}(z_1) \dots \Phi_{\ve_n}(z_n)
\end{equation}
has analytic matrix elements in the region $|z_1| \gg \dots \gg |z_n|$,
and extends to a meromorphic function in $(\C \setminus \{0\}).$ Its
highest component to highest component matrix elements satisfy the
quantum Knizhnik--Zamolodchikov
equation.
\end{thm}

\subsection{Construction of type I vertex operators.}
In \cite{JM} the type I vertex operator corresponding to the coproduct 
$\Delta$ is constructed.  We recall
\begin{prop}\cite{JM}
The vertex operators $\Phi_{\pm}$ have the following explicit
expressions in $\H$:
\begin{align}
& \label{phiminus}
\phioneminusi_{-}(z) = \expplusaminus{7}\expminusaplus{-5}e^{\alpha/2}
(-q^3z)^{(\partial + I/2)}, \\
& \label{phiplus}
\pP^{(1-i,i)} (z) = \pM^{(1-i,i)}(z)  x_0^- - q x_0^- \pM^{(1-i,i)}(z),
\end{align}
where $I$ is the operator on $V(\Lambda_i)$ such that $I(i) = i.$
\end{prop}

In order to calculate the bar action on $\H$, we define $\ov{\Delta}$ to
be the coproduct conjugated with the bar operator \eqref{bar}.
Explicitly:
\begin{equation}
\begin{split}
& \ov{\Delta}(K_i) = K_i \otimes K_i, \ \ov{\Delta}(E_i) = K_i^{-1}
\otimes E_i + E_i \otimes 1, \\ & \ov{\Delta}(F_i) = F_i \otimes K_i
+ 1 \otimes F_i, \ \ i = 0, \ 1.
\end{split}
\end{equation}
We now construct the vertex operators $\bpP,\ \bpM$ corresponding
to the coproduct $\ov{\Delta}.$

Let $\hat x_k^+ = x_k^+ K,\ \hat x_k^- = K^{-1} x_k^-$.
Using \eqref{interproperty} with $x$ equal to $K, E_i, F_i\ (i = 0, 1)$
respectively, the following must hold:
\begin{align}
& K \bppm K^{-1} = q^{\mp1}\bppm,  \label{kphi}\\
& \bpP(z) \hxzp = q \hxzp \bpP(z) + q^{-1}\bpM(z) , \label{phiplusxzeroplus} \\
& \bpM(z) \hxzp = q^{-1} \hxzp \bpM(z), \label{phiminusxzeroplus} \\
&\bpP(z) \hxzm = \hxzm \bpP(z), \label{phiplusxzerominus} \\
&\bpM(z) \hxzm = \hxzm \bpM(z) + q K^{-1} \bpP(z), \label{phiminusxzerominus}\\
&\bpP(z) \hxmp = q^{-1} \hxmp \bpP(z) + q^{-2} z^{-1}
\bpM(z), \label{phiplusxminusoneplus}  \\
&\bpM(z) \hxmp = q \hxmp \bpM(z), \label{phiminusxminusoneplus}\\
&\bpP(z) \hxpm = \hxpm \bpP(z), \label{phiplusxoneminus} \\
&\bpM(z) \hxpm = \hxpm \bpM(z) + z \bpP(z) K_1. \label{phiminusxoneminus}
\end{align}

We show these relations together with \eqref{normalization}
determine $\bppm(z)$ uniquely.  Assuming \er{kphi} through
\er{phiminusxoneminus} we show
\begin{lem}
\begin{align}
&\Xm{w} \bpP(z) - q \bpP(z) \Xm{w}  = 0, \label{xmbarphiplus} \\
&[a_n, \bpP(z)] = \frac{- q^{-n/2}}{q^n + q^{-n}} z^n \bpP(z),
\label{anphiplus} \\
&[a_{-n}, \bpP(z)] = \frac{- q^{-n/2}}{q^n + q^{-n}} z^{-n} \bpP(z)
\label{aminusnphiplus}.
\end{align}
\end{lem}

\begin{pf}
For the  proof of \eqref{anphiplus} and
\eqref{aminusnphiplus}, first calculate:
\begin{align*}
&[a_1,\bpP] = [\frac{C^{1/2} K^{-1}}{q+q^{-1}}[\xzp, \xom], \bpP] =
\fracqqi{C^{1/2}}(K^{-1}[\xzp,\xom]\bpP - \bpP K^{-1} [\xzp,\xom]) \\
&= \fracqqi{C^{1/2}K^{-1}} ([\xzp,\xom] \bpP - q^{-1} \bpP[\xzp,\xom])
\\
&= \fracqqi{C^{1/2}K^{-1}} (\xzp \xom \bpP - \xom \xzp \bpP - q^{-1}
\bpP \xzp \xom + q^{-1} \bpP \xom \xzp) \\
&= \fracqqi{C^{1/2}K^{-1}} (q^{-1}\xzp \bpP \xom - \xom(\bpP \xzp -
K^{-1} \bpM) - q^{-1} \bpP \xzp \xom + q^{-1} \bpP \xom \xzp) \\
&= \fracqqi{C^{1/2}K^{-1}} (-q^{-1}K^{-1} \bpM \xom + q^{-2} K^{-1} \xom
\bpM) \\
&=  \fracqqi{C^{1/2} K^{-2}} q^{-1} (- \bpM \xom + q^{-1} \xom \bpM)
= - \fracqqi{C^{1/2} q^{-1}} \bpP z = - \fracqqi{q^{-1/2} z} \bpP.
\end{align*}
Similarly we have:
\begin{equation*} [a_{-1}, \bpP] = \fracqqi{-q^{-1/2} z^{-1}}\bpP.
\end{equation*}
The relations \eqref{phiplusxzerominus} and \eqref{phiplusxoneminus}
together with repeated use of \er{aandx} now yield
\begin{equation} \Xm{w} \bpP(z) - q \bpP(z) \Xm{w}= 0.
\label{phiplusbarzxminusw}
\end{equation}
Let
\begin{align*}
&A(w) = \exp((q-q^{-1}) \sum_{n=1}^{\infty} \frac{[2n]}{n}
a_n q^{-n/2} w^{-n}), \\
&B(w) = \exp(- (q-q^{-1}) \sum_{n=1}^{\infty}\frac{[2n]}{n}
a_{-n} q^{-n/2} w^n).
\end{align*}
We have the following two identities:
\begin{align}
& (q-q^{-1}) [X^-(w),\xzp] = -w^{-1}KA(w) + w^{-1} K^{-1} B(w), \\
& (q-q^{-1}) [X^-(w),\xmp] = -w^{-2}(-q^{-1}KA(w) +  q K^{-1} B(w)).
\end{align}
Using \eqref{phiplusxzeroplus} and \eqref{phiplusxminusoneplus} we
obtain
\begin{equation}
qz(\bpP \hxmp - q^{-1} \hxmp \bpP) - (\bpP \hxzp - q \hxzp \bpP) = 0.
\label{prev}
\end{equation}
Bracketing \eqref{prev} with $(q-q^{-1}) X^-(w)$ we obtain the identity
\begin{equation}
\begin{split} \label{master}
& qz\bigl(\bpP(-q^{-1}KA(w) + qK^{-1}B(w))w^{-2} - q^{-1}w^{-2}(-q^{-1}K A(w)
+  \\ &q K^{-1} B(w))\bpP\bigr) -
 \bigl(\bpP (-w^{-1}KA(w) + w^{-1} K^{-1}B(w)) \\ & - q(-w^{-1} K A(w) + w^{-1}
K^{-1} B(w)) \bpP\bigr) = 0.
\end{split}
\end{equation}
Expanding this identity with respect to $w$ we obtain \er{anphiplus} and
\er{aminusnphiplus} to complete the proof.
\end{pf}
\begin{prop}
The commutation relations \er{xmbarphiplus}, \er{anphiplus},
\er{aminusnphiplus}, and \er{phiplusxzeroplus} determine the operators
$\bppm(z).$
\begin{align}
& \label{phiplusbar}
\bpP^{(1-i,i)}(z) = \expminusaminus{-}\expplusaplus{-} \\ & \times
e^{-\alpha/2}
(-qz)^{(-\partial + i)/2}(-q)^{i-1}, \\
& \label{phiminusbar} \bpM^{(1-i,i)}(z) = K [\bpP(z), x_0^+].
\end{align}
\end{prop}

\begin{pf}
The equations \eqref{master} and \eqref{phiplusbarzxminusw}, together
with the normalization \eqref{normalization} for $\bpM(z)$
determine the operator $\bpP(z)$
uniquely as given in \eqref{phiplusbar}.  To complete the proof, the
remaining relations must be checked.  For \eqref{phiminusxminusoneplus}
and \eqref{phiminusxzeroplus}, we use the following formulas, which are
derived from the defining expression for $\bpM$
\eqref{phiplusxzeroplus}.

\begin{align}
  & \label{norordXplusbpP} X^+(w) \bpP(z) = \frac{1}{w-q^{-1}z} \
  :\Xp{w} \bpP(z): , |z| \le |qw|, \\
& \label{norordbpPXplus}  \bpP(z) X^+(w) = \frac{1}{w-qz}\
  :\bpP(z) X^+(w):, |q^{-1}w| \le |z|, \\
& \label{norordXplusXplus}  \Xp{w_1} \Xp{w_2} = w_1^2
  (1-\frac{w_2}{w_1}) (1 - \frac{q^{-2}w_2}{w_1}) \ :\Xp{w_1} \Xp{w_2}:.
\end{align}
Here we have introduced the normally ordered
product on the $\{ a_n, \partial, \alpha | n \in \Z \setminus \{0\} \}$ by
\begin{align}
\label{normal}
:a_k a_l :\  &= a_k a_l \text{ if } k<0, \\
           &= a_l a_k \text{ if } k>0, \notag \\
 :\alpha \partial: \  &= \ : \partial \alpha: \  = \alpha \partial. \notag
\end{align}

We note the following identities, where the path of integration is the boundary
of suitable  two dimensional disk about the origin:
\begin{align*}
  & x_m^+ \bpP(z) x_n^+ = \frac{1}{(2\pi i)^2} \int \int dw_1 dw_2 w_1^m
  w_2^n \Xp{w_1} \bpP(z) \Xp{w_2} \\ & = \frac{1}{(2\pi i)^2} \int \int
  dw_1 dw_2 w_1^m w_2^n \frac{1}{w_1 - q^{-1}z} w_1^2 (1 -
  \frac{w_2}{w_1}) \\ & \hskip .8 in \times (1 - \frac{q^{-2}w_2}{w_1})
  \frac{1}{w_2 - qz} :\Xp{w_1} \Xp{w_2} \bpP(z): \\ & = \frac{1}{(2\pi
    i)^2} \int \int dw_1 dw_2 w_1^m w_2^n q^{-2} w_1 (1 -
  \frac{w_2}{w_1}) \\ & \hskip .8 in \times \bigl[\frac{q^2}{w_2 - qz} -
  \frac{1}{w_1 - q^{-1}z}\bigr] :\Xp{w_1} \Xp{w_2} \bpP(z): \\
  \intertext{Since exchanging $w_1$ and $w_2$ has no effect on
    $:\Xp{w_1} \Xp{w_2} \bpP(z):$ we make the variable substitution in
    the second summand to get:} & = \frac{1}{(2\pi i)^2} \int \int dw_1
  dw_2 \bigl[ \frac{w_2^m w_1^n q^{-2}w_2 (1-w_1/w_2) q^2}{w - qz} -
  \frac{w_1^m w_2^n q^{-2} w_1 (1 - w_2/w_1)}{(w_1 - q^{-1} z)} \bigl]\\
  & \hskip 1in \times :\Xp{w_1} \Xp{w_2} \bpP(z):.
\end{align*}
And similarly:
\begin{align*}
  & x_m^+ x_n^+ \bpP(z) = \frac{1}{(2\pi i)^2} \int \int dw_1 dw_2
  \Xp{w_1} \Xp{w_2} \bpP(z) w_1^m w_2^n \\ & = \frac{1}{(2\pi i)^2} \int
  \int dw_1 dw_2 \frac{q^{-2} w_1^m w_2^{n+1} - w_1^{m+1} w_2^n - q^{-2}
    w_1^{n+1} w_2^m + w_2^{m+1} w_1^n}{w_1 - q^{-1}z}. \\
  \intertext{And} & \bpP(z) x_m^+ x_n^+ = \frac{1}{(2\pi i)^2} \int \int
  dw_1 dw_2 \bpP(z) \Xp{w_1} \Xp{w_2} w_1^m w_2^n \\ & = \frac{1}{(2\pi
    i)^2} \int \int dw_1 dw_2 \frac{q^{-2} w_1^m w_2^{n+1} - w_1^{m+1}
    w_2^n - q^{-2} w_1^{n+1} w_2^m + w_2^{m+1} w_1^n}{w_1 - qz}.
\end{align*}

A formula for \eqref{phiminusbar} is determined
explicitly via \eqref{phiplusxzeroplus}.
Now \eqref{phiminusxminusoneplus} and
\eqref{phiminusxzeroplus} are seen to hold by
using the defining expression for $\bpM$
\eqref{phiplusxzeroplus}.  Also, we can now check explicitly
\eqref{phiminusxzerominus} and \eqref{phiminusxoneminus}. For example,
\begin{align*}  [\bpM(z), \hxpm] & = [q \bpP(z) \hxzp - q^2 \hxzp \bpP(z),
  \hxpm] \\ & = [2] C^{-1/2} \bigl( q \bpP(z) K a_1 - q^2 K a_1 \bpP(z) \bigr)
  = z \bpP(z) K,
\end{align*}
gives \eqref{phiminusxminusoneplus}.
\end{pf}

\subsection{Vertex operator action on the basic representation.}
We consider the action of $\ppm(z), \ \bppm(z)$ on $\H_\A = V(\Lzero)_\A
\oplus V(\Lone)_\A.$  
This is given as follows:
\begin{prop} \label{fouractions}
  Let $m \ge 0, i = 0,1.$
\begin{align}
& (-q)^{\partial-1} \bpM(z)  \eali{m} = \expplusaminus{-5} \expminusaplus{3}
 \label{bpMaction}  \\ & \hskip .8 in
\times
e^{\alpha/2}   (-q)^{(\partial + 3I)/2}
z^{(\partial + I)/2}   \eali{m},  \text{where if $m=0$, $i=1$,} \notag\\
& (-q)^{1-\partial} \pM(z)   \eali{\pm m} =
\expplusaminus{7} \expminusaplus{-5}
\label{pMaction} \\ & \hskip .8 in
\times
e^{\alpha/2} (-q)^{(\partial + 3I)/2}
z^{(\partial + I)/2}  \eali{\pm m}. \notag
\end{align}
The action of $\bpP$ and $\pP$
is as follows:
\begin{align}
& (-q)^{-\partial} \bpP(z)  \eali{\pm m}= \expminusaminus{-}
\expplusaplus{-} \label{bpPaction} \\ & \hskip .8 in
\times
e^{-\alpha/2}  (-q)^{(-3\partial + 3I)/2}
z^{(-\partial + I)/2}  \eali{\pm m},  \notag \\
& (-q)^\partial \pP(z)  \eali{-m}= \expminusaminus{11} \expplusaplus{-9}
 \label{pPaction} \\ & \hskip .8 in
\times
e^{-\alpha/2}  (-q)^{(-3\partial + 3I)/2}
z^{(-\partial + I)/2}  \eali{-m}, \text{where if $m=0$, $i=0$}. \notag
\end{align}
\end{prop}

\begin{pf} The proof of \er{pMaction} follows directly from
\er{phiminus}.  We consider \er{bpMaction}.
\begin{equation*}
\begin{split} \bpM(z) & e^{m \alpha} \vi  = K [\bpP(z), \xzp] e^{m \alpha}
  \vi = K \frac{1}{2\pi i} \int dw [\bpP(z), \Xp{w}] e^{m \alpha} \vi 
  \\ \intertext{using \er{norordXplusbpP} and \er{norordbpPXplus}} & = - K
  \frac{1}{2\pi i} \int \frac{dw}{(w-q^{-1}z)} \ :\bpP(z) \Xp{w}: e^{m
    \alpha} \vi  \\ & =  -K :\bpP(z) \Xp{q^{-1}z}: e^{m \alpha} \vi  \\
  & =  \expplusaminus{-5}\expminusaplus{3} \\ & \hskip .1in 
   \times e^{\alpha/2} (-q)^{(3I
       -   \partial)/2} z^{(\partial + I)/2} \eali{m}, 
\end{split}
\end{equation*}
which gives \er{bpMaction}. \er{bpPaction} and \er{pPaction} are similar.
\end{pf}


\section{Macdonald Polynomials}

\subsection{Partitions}
As usual, by a {\em partition} we mean a sequence of non--negative
integers in decreasing order
\begin{equation} \lambda_1 \ge \lambda_2 \ge \dots \ge \lambda_r \ge
  \dots
\end{equation}
containing finitely many non--zero terms.
The number of non--zero $\lambda_i$ is called the {\em length} of $\lambda$,
denoted by $\ell(\lambda)$, and each $\lambda_i$ is called a {\em part}.
The sum $|\lambda | \  =
\sum_i\lambda_i$ is called the {\em weight} of
$\lambda.$   Given a partition $\lambda$, if
the part $i>0$ appears $m_i$ times we write
\begin{equation*} \lambda = (1^{m_1}2^{m_2} \dots r^{m_r} \dots).
\end{equation*}
We define the integer $$z_\lambda = \prod_{i \ge 1} i^{m_i} m_i!.$$
The {\em dual partition} $\lambda'$ is defined by  setting its parts as
\begin{equation*} \lambda'_i =\ \text{Card}\{j: \lambda_j \ge i\}.
\end{equation*}
We define the {\em dominance partial ordering}
on partitions by setting $\lambda \le \mu$
if $\sum_i \lambda_i = \sum_i \mu_i$ and
$\lambda_1 + \dots + \lambda_k \le \mu_1 + \dots + \mu_k$ for every $k.$
That this is a partial order is made clear by considering the two partitions
$(1^3,3), \ (2^3),$  which are incomparable.

Denote by $\symfun$ the ring of symmetric functions in countably
many variables $\{x_i \ | \ i \ge 1\}$.  Let $S^\infty$ be the
permutations of $\N$ which fix a cofinite set.
For a partition $\lambda$  define the {\em monomial symmetric functions}
\begin{equation} m_\lambda = \sum_{\alpha = s \lambda, \ s \in S^\infty}
  \prod_{i} x_i^{s(\lambda_i)}.
\end{equation}

For each $r \ge 1$ the $r$--th power sum is
\begin{equation} p_r = \sum_{i \ge 1} x_i^r.  \label{powersum}
\end{equation}
The $p_r \  (r \ge 1)$ 
form a polynomial basis of $\symfun_\Q = \symfun \otimes \Q$.
Define
\begin{equation} \label{plambda} p_\lambda = p_{\lambda_1} p_{\lambda_2} \dots
\end{equation}
for each partition $\lambda = (\lambda_1, \lambda_2, \dots).$
The $p_\lambda, \ \lambda$ a partition, form a $\Q$--basis of
$\lambda^{\Sigma}_{\Q}.$  
\subsection{Properties of Macdonald polynomials}

We recall some basic facts about Macdonald polynomials.  We refer the
reader to Chapter VI of \cite{Ma} for further information.
Let $\C(q,t)$ be the field of rational functions in $q$ and $t$.
Define a scalar product $( \cdot \ , \ \cdot )$ on
$\symfun_\Q$ by
\begin{equation}
( p_\lambda, p_\mu ) =  ( p_\lambda, p_\mu )_{q,t} =
\delta_{\lambda,\mu} z_\lambda \prod_{i = 1}^{l(\lambda)} \frac{1 -
  q^{\lambda_i}}{1 - t^{\lambda_i}}.
\end{equation}
Let $\{x_i\}_{i \ge 0}$ be an infinite set of indeterminates.
\begin{thm} \cite{Ma} Let $\lambda$ be a partition. There exists a unique
family of symmetric functions
$P_\lambda(x;q,t) \in \C(q,t)[x_1, x_2, \dots]$ which satisfy the
following properties:
\begin{enumerate}
\item $P_\lambda$ is symmetric with respect to the $x_i$, $i \ge 1$.
\item $P_\lambda = m_\lambda + \sum_{\lambda < \mu} K_{\lambda, \mu}
  m_\mu,$ where $K_{\lambda, \mu} \in \Q(q,t).$
\item The $P_\lambda$ are pairwise orthogonal relative to the scalar
  product $( \cdot \ | \ \cdot )$ and $( P_\lambda,
  P_\lambda ) = b^{-1}_\lambda(q,t).$
\end{enumerate}
\end{thm}
Here
\begin{align} \label{blambda}
b_\lambda(q,t) = \prod_{s \in \lambda} \frac{1 - q^{a(s) +
   1}t^{l(s)}}{1 - q^{a(s)}t^{l(s) + 1}},
\end{align}
 where for each square $s = (i,j)$ in the diagram of $\lambda$ we have
\begin{align*}
  & a(s) = a_\lambda(s) = \lambda_i - j, \\
  & l(s) = l_\lambda(s) = \lambda'_j - i.
\end{align*}
We also denote the numerator and denominator of $b_\lambda$ by
$c_\lambda(q,t)$ and $c'_\lambda(q,t)$ respectively.  
We refer to the form above as {\em Macdonald's form}.  Setting 
$Q_\lambda(x;q,t) = b_\lambda(q,t) P_\lambda(x;q,t)$ we see $Q_\lambda$
is dual to $P_\lambda$ with respect to Macdonald's form. 
If $z$ is an indeterminate define
\begin{equation*} (z;q)_\infty = \prod_{j \ge 0} (1 - z q^j), \quad
 \xi(z) = \frac{(q^2 z; q^4)_\infty}{(q^4 z;
    q^4)_\infty}.
\end{equation*}
A direct calculation shows:
\begin{equation*} \exp\bigl(\sum_{n = 1}^\infty \frac1{n}
  \frac{[n]}{[2n]} z^n \bigr)
= \frac{(q^3z;q^4)_\infty}{(qz;q^4)_\infty}.
\end{equation*}
and from here we have immediately
\begin{equation} \label{baridentity}
\ov{\xi(z)} = (1-z)/\xi(z), 
\end{equation}
where $\ov{\ }$ is extended to rational functions in $z$ by setting
$\ov{qz} = q^{-1}z.$ 
It is known that the $P_\lambda$ can be expressed via generating
series as follows:
\begin{align} \label{innerproductPi} \Pi(x,y;q,t) & :=
\prod_{i,j} \frac
{(tx_iy_j;q)_\infty}{(x_iy_j;q)_\infty} = \exp\bigl(\sum_{n \ge 0} 
   \frac{1-t^n}{1-q^n} \frac{1}{n}
\sum_i x_i^n \sum_j y_j^n\bigr) \\ & = \sum_\lambda P_\lambda(x_i;q,t)
Q_\lambda(y_j;q,t). 
\end{align}
The Macdonald polynomials also satisfy: 
\begin{align} & P_\lambda(x;q,t) = P_\lambda(x;q^{-1},t^{-1}),
  \label{inverseproperty} \\
              &  Q_\lambda(x;q,t) = (qt^{-1})^{|\lambda|} 
Q_\lambda(x;q^{-1},t^{-1}), 
\end{align}

The power sums $p_r = \sum x_i^r$ form a $\Q$--basis of the ring of
symmetric functions in the $x_i$.
\begin{defn} Let $\P(p_n;q,t)$  to be the  $\Q(q,t)$
polynomial in $p_i,\  i >0$ for which $\P(p_n;q,t) = P(x_n;q,t).$
\end{defn}
\begin{prop} Let $N>0$. Let $z_i,\ i \ge 0$ be indeterminates.
The Macdonald polynomials satisfy the
following generating series:
\begin{equation} \label{generatingfunction}
        \exp\bigl(\sum_{n\ge1}(-1)^{n-1} \frac{b_n}{n}
  \sum_{j=1}^N z_j^n \bigr) = \sum_{l(\lambda) \le N} P_\lambda(z;q,t)
\P_{\lambda'}(b_{n};t,q).
\end{equation}
\end{prop}
\begin{pf}
This follows from \cite{Ma} (page 310), where we restrict to the
ring of symmetric polynomials in $z_1, \dots, z_N.$
\end{pf}

Finally, introduce the involution of $\C(q,t)[x_1, x_2, \dots]$ by
setting
\begin{equation}
\omega_{q,t} (p_n) = \frac{1-q^n}{1-t^n} (-1)^{n-1} p_n.
\end{equation}
Then (c.f. \cite{Ma} VI 5.1), 
\begin{equation} \omega_{q,t} P_{\lambda}(x;q,t) = Q_{\lambda'}(x;t,q).
\end{equation}


\section{Invariance of Macdonald polynomials under bar action}
Following \eqref{rho} we induce a bar
action on $\H$ by
\begin{defn} Let $u \in \uqatwoatone$, for the highest weight vector $\vi$
of $V(\Lambda_i), i=0, 1,$
define $\ov{u \vi} = \ov{u} \vi.$
\end{defn}
This implicitly defines a bar action on $V(\Lambda_i),$ and it is clear that
\begin{prop}
For $\Phi(z): \H \rightarrow \H \otimes V(z)$ we have
\begin{equation} \ov{\Phi(z) u\vzero + v\vone} = \ov{\Phi}(z) (\ov{u} \vzero +
  \ov{v}\vone),
\label{bardefined}
\end{equation}
\end{prop}
where $u, v \in \uqatwoatone$.

In Proposition \ref{fouractions}, the application of a vertex operator
involves a multiplication by a power of $q$.  In order to simplify the
statement of results we introduce the following normalization:
\begin{defn}Let $m \ge 0, i
= 0, 1.$
\begin{align*}
&v_{m,i}  = (e^{\alpha/2} (-q)^{(\partial + 3I)/2})^{2m} \vi  =
(-q)^{m (m +1+i)} e^{m\alpha} \vi, \\
&v_{-m,i} =  (e^{-\alpha/2} (-q)^{(-3\partial + 3I)/2})^{2m} \vi  =
 (-q)^{3 m (m-i)}  e^{-m \alpha} \vi.
\end{align*}
where when $m = 0,$ we have $i=1$ in the first case and $i=0$ in the
second case.
\end{defn}
As usual, denote by $\lfloor n \rfloor$
 the largest integer less than or equal to
$n$.
\begin{prop}  \label{genseries}
Let $i=0,1, N \ge 1,\ n = \lfloor\frac{N+1-i}2 \rfloor,
\ m = \lfloor\frac{N+i}2\rfloor , k = \frac{1 + (-1)^{N+i+1}}2.$
\begin{align}
\begin{split}
 & \label{bpMgenseries}
(-q)^{\partial-1} \bpM(z_1)  (-q)^{\partial-1} \bpM(z_2) ...
 (-q)^{\partial-1} \bpM(z_N) \vi, \\
 & =
\frac{\prod_{i < j} (1 - z_j/z_{i})}{\prodxi} \prodz \expplusaminusz{-5}
 v_{m,k}, \\ & \text{ where if $m=0$, $i=1$,}
\end{split}  \\
\begin{split}
 &\label{pMgenseries}
 (-q)^{1-\partial} \pM(z_1) (-q)^{1-\partial} \pM(z_2) ...
 (-q)^{1-\partial} \pM(z_N) \vi \\ & =
{\prodxi} \prodz \expplusaminusz{7}  v_{m,k},
\end{split} \\
\intertext{and}
\begin{split}
  &  \label{bpPgenseries}
(-q)^{-\partial}  \bpP(z_1)  (-q)^{-\partial} \bpP(z_2) ...
 (-q)^{-\partial} \bpP(z_N) \vi
  =  \\ &   \prodzminus
\frac{\prod_{i < j} (1 - z_j/z_{i})}{\prodxi}  \expminusaminusz{-}
v_{-n,k}, \\ & \text{ where if $m=0$, $i=0$},
\end{split} \\
\begin{split}
 &\label{pPgenseries}
 (-q)^{\partial} \pP(z_1) (-q)^{\partial} \pP(z_2) ...
 (-q)^{\partial} \pP(z_N) \vi =
 \\ &  \hskip .8 in
\prodzminus {\prodxi} \expminusaminusz{11} v_{-n,k}.
\end{split}
\end{align}
\end{prop}
\begin{pf} 
These follow by direct calculation from Proposition \ref{fouractions}.
 \end{pf}
Let
\begin{equation*}
 c(z) =
{\prodxi} \prodz,  \quad
 d(z) =
{\prodxi} \prodzminus.
\end{equation*}

Using the generating function \er{generatingfunction}, the formulas
\er{bpMaction}, \er{pMaction},
\er{bpPaction} and \er{pPaction} give the following identities:
\begin{prop} \label{firstmac} Let $N, n, m, k$
 be  as in the previous proposition.
\begin{align}
\begin{split} \label{bpMmac}
 (-q)^{\partial-1}& \bpM(z_1)  (-q)^{\partial-1} \bpM(z_2) ...
 (-q)^{\partial-1} \bpM(z_N) \vi  \\ & =
\ov{c(z)} \macsum{\{q^{-3+(1/2)} z_j\}_{j =
     1}^N}{(-1)^{n-1}  a_{-n} }{4}{2} v_{m,k},
\end{split} \\
\begin{split}
 (-q)^{1-\partial}& \pM(z_1) (-q)^{1-\partial} \pM(z_2) ...
 (-q)^{1-\partial} \pM(z_N) \vi \\ & = c(z) \macsum{\{q^{3+(1/2)} z_j\}_{j =
     1}^N}{(-1)^{n-1}  a_{-n}}{4}{2} v_{m,k},
\end{split} \notag \\
\intertext{and}
\begin{split}
 (-q)^{-\partial} & \bpP(z_1)  (-q)^{-\partial} \bpP(z_2) ...
 (-q)^{-\partial} \bpP(z_N) \vi \\
 & = \ov{d(z)} \macsum{\{q^{-3 + (5/2)} z_j\}_{j =
     1}^N}{(-1)^{n} a_{-n}}{4}{2} v_{-m,k},
\end{split} \notag \\
\begin{split} \label{pPmac}
 (-q)^{\partial}& \pP(z_1) (-q)^{\partial} \pP(z_2) ...
 (-q)^{\partial} \pP(z_N) \vi \\
  & = d(z)  \macsum{\{q^{3+(5/2)}z_j\}_{j =
     1}^N}{(-1)^{n} a_{-n}}{4}{2} v_{-m,k}.
\end{split}
\end{align}
\end{prop}
From these and \ref{inverseproperty} we immediately obtain:
\begin{prop} \label{invariance}  Let $m \ge 0,\ i= 0, 1.$
Let $\Delta = 3 \min(2m + i - (\ell(\lambda') - 1), 0).$
\begin{align}
  &\ov{(q^{-\Delta})q^{|\lambda|/2}\P_{\lambda}( \ahat_{-n}; q^2, q^4)
    v_{m,i}} = (q^{-\Delta}) q^{|\lambda|/2}\P_{\lambda}( \ahat_{-n}; q^2,
  q^4)  v_{m,i}, \label{biplusm}\\
  &\ov{ (q^{\Delta}) q^{5|\lambda|/2}\P_{\lambda}( \ahatn_{-n};
    q^2, q^4) v_{-m,i} } = (q^{\Delta}) q^{5|\lambda|/2}
  \P_{\lambda}( \ahatn_{-n}; q^2, q^4) v_{-m,i}. \label{biminusm}
\end{align}
where when $m = 0$ we have $i=1$ in the first case and $i=0$ in the
second case.
\end{prop}
\begin{pf}
In both cases the bar--invariance follows directly from
Proposition \ref{firstmac} when $l (\lambda') \le \lfloor m-1/2\rfloor.$
An extra factor of $q^\Delta$ appears for $\mathcal
P_\lambda$ in sectors which don't satisfy this inequality.  Fix an
arbitrary partition $\lambda.$ By Proposition \ref{firstmac} for large
enough $m$, we have $\Delta = 0$.  Now apply the $\bpP(w)$ action of
\er{bpPaction} to both sides.  Using \er{bpMgenseries} we have
\begin{equation*}
\begin{split}
 &\bpP(w) (-q)^{\partial-1} \bpM(z_1)  (-q)^{\partial-1} \bpM(z_2) ...
 (-q)^{\partial-1} \bpM(z_N) \vi \\
 & =  \prod_{j=1}^N \frac{(z_j/w,q^4)_\infty}{(q^{-2}z_j/w,q^4)_\infty}
w^{-2m}
\frac{\prod_{i < j} (1 - z_j/z_{i})}{\prodxi} \prodz
\\
& \hskip .8 in \exp\bigl(-\sum_{n=1}^\infty \frac{q^{-n/2}}{n} a_{-n}
w^n\bigr)
\expplusaminusz{-5}
 v_{m-k,1-k}.
\end{split}
\end{equation*}
and from \er{pMgenseries} we have
\begin{equation*}
\begin{split}
 & \bpP(w) (-q)^{1-\partial} \pM(z_1) (-q)^{1-\partial} \pM(z_2) ...
 (-q)^{1-\partial} \pM(z_N) \vi \\ & =
\prod_{j=1}^N \frac{(q^6z_j/w,q^4)_\infty}{(q^{4}z_j/w,q^4)_\infty}
w^{-2m}
{\prodxi} \prodz
 \\
& \hskip .8 in
\exp\bigl(-\sum_{n=1}^\infty \frac{q^{-n/2}}{n} a_{-n} w^n\bigr)
\expplusaminusz{7} v_{m-k,1-k},
\end{split}
\end{equation*}
Noting that  \eqref{baridentity} implies
\begin{equation*}
\ov{\prod_{j=1}^N \frac{(z_j/w,q^4)_\infty}{(q^{-2}z_j/w,q^4)_\infty}}
= \prod_{j=1}^N \frac{(q^6z_j/w,q^4)_\infty}{(q^{4}z_j/w,q^4)_\infty}
\end{equation*}
we see that
\er{biplusm} follows for $\P_{\lambda}(\ahat_{-n}; q^2, q^4)
v_{m-k,1-k}$ by considering the
coefficient of $w^{-2m}$. A power of $q^3$ appears once for each $\bpP$
applied, which accounts for $q^\Delta$.  This implies \er{biplusm} for all $m$.
The proof of \er{biminusm} is similar using
$\pM(w)$ and \er{bpPgenseries} and \er{pPgenseries}.
\end{pf}


\section{The coincidence of Kashiwara's and Macdonald's forms}

In the following we show that the forms characterizing  the canonical
basis and the Macdonald polynomials coincide in our realization 
of the ring of symmetric functions in the representation 
$\H$.  As with the bar invariance in the previous section, the
$\mathcal P_\lambda$ are multiplied by a power of $q$, 
depending on the sign of $m$  sector determined by $v_{m,i}.$

 We note that the bilinear form of Proposition 1  on the representations
$V(\Lone)$ and $V(\Lzero)$ extends naturally to $\H$ by requiring
$V(\Lone)$ and $V(\Lzero)$ to be orthogonal.  Furthermore, there is also
a form with the property $\eqref{rhochange}$ on $V_z$.  This is obtained
by starting with such a form for the two dimensional representation
$V$ of $\uqtwo$ and extending it to $V_z = V_1 \otimes C[z,z^{-1}]$
by setting $(v \otimes z^n, w \otimes z^m) = (v,w)
\delta_{n,m}.$ 
\begin{defn} The form $(\ ,\ ): \H \otimes  V_z \times \H \otimes
 V_z \rightarrow \Q(q)$
is the unique bilinear form determined by $(x \otimes v z^n, y
 \otimes w z^m) = (x,y)_{\H} (v,w) \delta_{m,n}.$
\end{defn}
Explicit expressions for $\Phi_{\ve}^*(z)$ (see
\eqref{phioneminusistar})  in terms of
$\Phi_{\ve}(z)$ are readily calculated.  We recall:
\begin{prop} \label{phidualtophi} \cite{JM}
 Let $u,  v \in \H$.  We have
\begin{equation}
\begin{split}   \Phi_\ve^{*}(z) =
  (-q)^{i+(\ve-1)/2} \Phi_{-\ve}(q^{-2}z), \notag \\
      \ov{\Phi}_\ve^{*}(z) =  (-q)^{-i+(1-\ve)/2}
     \ov{\Phi}_{-\ve}(q^{2}z).
\end{split}
\end{equation}
\end{prop}
We have
\begin{prop} \label{ipswitch}
 Let $u,  v \in \H$.  We have
\begin{equation} \label{innerproductidentity}
\begin{split} (u,\Phi_\ve(z)v) = (\Phi_\ve^{*}(z^{-1})u, v),  \notag\\
(u,\ov{\Phi}_\ve(z)v) = (\ov{\Phi}_\ve^{*}(z^{-1})u, v). \notag
\end{split}
\end{equation}
\end{prop}
\begin{pf}
 A direct calculation shows that $\rho$ commutes with both coproducts:
 i.e. $\Delta \circ \rho = \rho \otimes
  \rho \circ \Delta$ and $\ov{\Delta} \circ \rho = \rho \otimes
  \rho \circ \ov{\Delta}$.  By definition, the
 vertex operators $\Phi$, $\Phi^*$ (see \eqref{phioneminusi},
 \eqref{phioneminusistar}) 
 are intertwiners for  $\Delta, \ov{\Delta}$ respectively.  Now we prove the
proposition by induction on the component degree of $\H \otimes V_z$.
For $u,v$ equal to $v_0$ or $v_1$ the statement easily verifiable.
We check that if the
proposition holds for $v \in \H$ with arbitrary $u \in \H \otimes V_z$,
then it holds for $a v$ where $a \in \uqatwo.$
Let $u_1, u_2 \in \H$, $v \in V_z.$
\begin{align*}
  (u_1 \otimes v, \Phi(z) a u_2) & = (u_1 \otimes v,
  \Delta(a)\Phi(z)u_2) = ((\rho \otimes \rho) \Delta(a) (u_1 \otimes v),
  \Phi(z) u_2) \\ & = (\Delta \rho(a) (u_1 \otimes v), \Phi(z) u_2) =
  (\Phi^*(z^{-1})(\Delta \rho(a)) (u_1 \otimes v), u_2) = \\ & =
  (\rho(a) \Phi^*(z^{-1}) (u_1 \otimes v), u_2) = (\Phi^*(z^{-1})u_1
  \otimes v, au_2).
\end{align*}
We note that when $u_2$ is $\Lone$ or $\Lzero$
the proposition clearly holds, and this completes the induction.
\end{pf}

\begin{lem} \label{iplemma} Let $m \ge 0, i = 0,1.$
\begin{align}
 &  {\bpM}^{*}(z^{-1}) (-q)^{2\partial-2}
\bpM(w) v_{m,i} = \frac{(q^{-8}zw,q^{-4})_\infty}
{(q^{-6}zw,q^{-4})_\infty} \tag{a} \\ & \hskip .8 in  \times
 \expminusaminuszsame{3}
\exp\bigl(-\sum_{n=1}^\infty \frac{q^{11n/2}}{n} a_{-n} w^n \bigr)
      v_{m,i}, \notag \\
 & {\pP}^{*}(z^{-1}) (-q)^{2\partial}
\pP(w) v_{-m,i} = \frac{(q^{8}zw,q^{4})_\infty}{(q^{6}zw,q^{4})_\infty}
\tag{b} \\ &
\hskip .8
in   \times
\exp\bigl(\sum_{n=1}^\infty \frac{q^{3n/2}}{n} a_{-n} z^{-n} \bigr)
\exp\bigl(\sum_{n=1}^\infty \frac{q^{-5n/2}}{n} a_{-n} w^n \bigr)
v_{-m,i}. \notag
\end{align}
\end{lem}
\begin{pf} This is a direct calculation using Proposition
  \ref{fouractions}.
\end{pf}
\begin{prop} \label{innerprod}
 Let $m \ge 0, i = 0,1.$
Let $\Delta = 3 \min(2m + i - (\ell(\lambda') - 1), 0).$ Then
\begin{align}
 & \bigl(q^{|\lambda|/2-\Delta} \P_{\lambda'}((-1)^{n-1}a_{-n};
  q^{-2},q^{-4}) v_{m,i} , q^{|\mu|/2 - \Delta}
             \P_{\mu'}((-1)^{n-1}a_{-n}; q^{-2},q^{-4})v_{m,i}\bigr)
             \tag{a} \\ &
             \hskip 1 in
  =\delta_{\lambda,\mu}  b_\lambda(q^{-4},q^{-2}) =  \delta_{\lambda,\mu}
b^{-1}_{\lambda'}(q^{-2},q^{-4}), \notag \\
 & \bigl(q^{5|\lambda|/2+ \Delta} \P_{\lambda'}((-1)^{n}a_{-n};
  q^{2},q^{4})v_{-m,i} ,
         q^{5|\mu|/2 + \Delta}\P_{\mu'}((-1)^{n}a_{-n}; q^{2},q^{4})v_{-m,i}
         \bigr) \tag{b}  \\ & \hskip 1 in
         =\delta_{\lambda,\mu}  b_\lambda(q^{4},q^{2}) = \delta_{\lambda,\mu}
 b^{-1}_{\lambda'}(q^{2},q^{4}).  \notag
\end{align}
where when $m = 0$ we have $i=1$ in the first case and $i=0$ in the
second case.
\end{prop}
\begin{pf}  We check (b).  As in the previous section  first 
we will restrict to the case where $l(\lambda) < \lfloor (m-i)/2
\rfloor$ and then we will extend to the general case.
\begin{align*}
  & \bigl( (-q)^{\partial}  \pP(z_1) (-q)^{\partial} \pP(z_2) ...
 (-q)^{\partial} \pP(z_N) \vi,  \\ & \hskip .8in (-q)^{\partial} \pP(w_1)
(-q)^{\partial} \pP(w_2) ... (-q)^{\partial} \pP(w_N) \vi  \bigr) \\
  & = \bigl( d(z)  \macsum{\{q^{3+(5/2)}z_j\}_{j =
     1}^N}{(-1)^{n} a_{-n}}{4}{2} v_{-m,j}, \\ & \hskip .8in d(w)
\macsummu{\{q^{3+(5/2)}w_j\}_{j =
     1}^N}{(-1)^{n} a_{-n}}{4}{2} v_{-m,j} \bigr) \\
   & =  d(z) d(w) \sum_{\lambda,\mu} q^{5/2(|\mu| + |\lambda|)}
P_\lambda(q^3z_i;q^4,q^2) P_\mu(q^3w_j;q^4,q^2) \\ & \hskip .8in \times
\bigl( \P_{\lambda'}((-1)^{n} a_{-n};q^2,q^4) v_{-m,j},
\P_{\mu'}((-1)^{n} a_{-n};q^2,q^4) v_{-m,j}  \bigr),
\end{align*}
by  \eqref{pPmac}.
On the other hand, by Lemma \ref{iplemma} we have,
\begin{align*}
  & = \bigl( (-q)^{\partial} \pP(z_1) (-q)^{\partial} \pP(z_2) ...
 (-q)^{\partial} \pP(z_N) \vi,  \\ & \hskip .8in (-q)^{\partial} \pP(w_1)
(-q)^{\partial} \pP(w_2) ... (-q)^{\partial} \pP(w_N) \vi  \bigr) \\
& = d(z) d(w) \prod_{1 \le i,j \le
  n}\frac{(q^8z_iw_j;q^4)_\infty}{(q^6z_iw_j;q^4)_\infty} (\vi, \vi) \\
& =   d(z) d(w) \macsumq{q^3z_i}{q^3w_j}{4}{2}.
\end{align*}

Now we check the general case.  Fix a partition $\lambda$.
We know that for $m$ large enough part
(b) holds.  Pick the largest $m$ for which
the result doesn't hold for $\P_{\lambda}$
and $v_{-m+1,j},$ where $j = 0,1.$ Applying Lemma \ref{iplemma}, we have
\begin{align*}
 &  \bigl( (-q)^{\partial} \pM(z_0)
(-q)^{\partial} \pP(z_1) (-q)^{\partial} \pP(z_2) ...
 (-q)^{\partial} \pP(z_N) \vi,  \\ & \hskip .8in (-q)^{\partial}
 \pM(w_0) (-q)^{\partial}  \pP(w_1)
(-q)^{\partial} \pP(w_2) ... (-q)^{\partial} \pP(w_N) \vi  \bigr) \\
  & = \bigl( {\tilde d}(z)
       \Nplusmacsum{q^{3+(5/2)}{\tilde z_0,\{q^{3+(5/2)}z_j\}_{j =
     1}^N}}{(-1)^{n} a_{-n}}{2}{4} v_{-m+1,j}, \\ & \hskip .2in {\tilde d}(w)
\Nplusmacsummu{q^{3+(5/2)}\tilde w_0,\{q^{3+(5/2)}w_j\}_{j =
     1}^N}{(-1)^{n} a_{-n}}{2}{4} v_{-m+1,j} \bigr) \\
   & =  {\tilde d}(z) {\tilde d}(w)
  \sum_{\lambda,\mu} q^{5/2(|\mu| + |\lambda|)}
P_\lambda(q^3 {\tilde z_0},q^3z_i;q^4,q^2)
         P_\mu(q^3 {\tilde w_0},q^3w_j;q^4,q^2)
\\ & \hskip .8in \times
\bigl( \P_{\lambda'}((-1)^{n} a_{-n};q^2,q^4) v_{-m+1,j},
\P_{\mu'}((-1)^{n} a_{-n};q^2,q^4) v_{-m+1,j}  \bigr), \end{align*}
where $\tilde w_0 = q^{-2} w_0$ and $\tilde z_0 = q^{-2} z_0$ and
\begin{equation*} {\tilde d}(z) = d(z) \prod_{i=1}^N
 \frac{(q^6 z_i/z_0; q^4)}{(q^4 z_i/z_0; q^4)}.
\end{equation*}
But, as before, the left hand side also equals:
\begin{align*}
  & \bigl( (-q)^{\partial} \pM(z_0)
(-q)^{\partial} \pP(z_1) (-q)^{\partial} \pP(z_2) ...
 (-q)^{\partial} \pP(z_N) \vi,  \\ & \hskip .6in (-q)^{\partial}
 \pM(w_0) (-q)^{\partial}  \pP(w_1)
(-q)^{\partial} \pP(w_2) ... (-q)^{\partial} \pP(w_N) \vi  \bigr) \\
  & = \bigl( (-q)^{\partial} \pP(z_1) (-q)^{\partial} \pP(z_2) ...
 (-q)^{\partial} \pP(z_N) \vi,  \\ & \hskip .6in
\pM^*(z_0^{-1})
(-q)^{2\partial}
 \pM(w_0) (-q)^{\partial}  \pP(w_1)
(-q)^{\partial} \pP(w_2) ... (-q)^{\partial} \pP(w_N) \vi  \bigr) \\
  & = {\tilde d}(z) {\tilde d}(w) \frac{(q^8 {\tilde z}_0 {\tilde w}_0;
    q^4)_\infty}{(q^6 {\tilde z}_0 {\tilde w}_0;q^4)_\infty} \times
        \prod_{1 \le i,j \le
  n}\frac{(q^8z_iw_j;q^4)_\infty}{(q^6z_iw_j;q^4)_\infty} (\vi, \vi) \\
  & = \tilde d(z) \tilde d(w)
             \Nplusmacsumq{q^3 {\tilde z}_0, q^3z_i}{q^3 {\tilde
      w}_0, q^3w_j}{4}{2}.
\end{align*}
The result now follows as above.
\end{pf}


\section{Lattice of dual Macdonald Polynomials}

Let $\hat a_{-n} = \omega_{q^4,q^2} (a_{-n}) = (-q)^n \frac{[2n]}{[n]}
a_{-n}.$   Then, changing variables, we have:
\begin{equation} \P_{\lambda'}(a_{-n};q^2,q^4) = {\mathcal Q}_{\lambda}(\hat
  a_{-n}; q^4, q^2).
\end{equation}
\begin{defn}  Let
\begin{align}
\tilde {\mathcal Q}_\lambda v_{m,i} =
\begin{cases} q^{|\lambda|/2} {\mathcal Q}_{\lambda'}(\hat a_{-n}; q^{4},q^{2})
v_{m,i} & \text{if $m>0$ or $m=0$ and $i=1$}, \\
   q^{5 |\lambda|/2} {\mathcal Q}_{\lambda'}(\hat a_{-n}; q^{-4},q^{-2})
v_{m,i} &\text{if $m<0$ or $m=0$ and $i=0$}.
\end{cases}
\end{align}
\end{defn}

As a result of the previous section we have:

\begin{prop} \label{orthogonal} The basis of the
  $\tilde {\mathcal Q}_\lambda v_{m,i}$ is orthogonal with respect to
   $(\ .\ ,\ .\ )$ and has dual $ b_\lambda(q^4,q^2) \tilde \P_\lambda
  v_{m,i} = \tilde {\mathcal Q}_\lambda v_{m,i}.$
\end{prop}

\begin{defn} Let $J^*_\lambda(q,t) = (c'_\lambda(q,t))^{-1} P_\lambda(q,t)
  = (c_\lambda(q,t))^{-1} Q_\lambda(q,t)$ be the dual integral
  Macdonald polynomials.  Let $\tilde {\mathcal J}$ (resp. $\tilde {\mathcal
    J^*}$) be $c_\lambda(q^4,q^2) \tilde \P$ (resp. $c'_\lambda(q^4,q^2)^{-1}
  \tilde \P$).
\end{defn}
Then we have
\begin{equation*}
(J_\lambda(q,t), J^*_\mu(q,t)) = \delta_{\lambda,\mu}.
\end{equation*}
where $J_{\lambda}(q, t)$ is the integral Macdonald function
(\cite{Ma}, p. 352) defined by
$J_{\lambda}(q, t)=c_{\lambda}(q, t)P_{\lambda}(q, t)=c_{\lambda}'(q, t)
Q_{\lambda}(q, t)$.
Let
\begin{equation}
{\mathcal L} = \oplus_{\A}    \tilde
{\mathcal J}^*_\lambda(a_{-n};q^2,q^4)v_{m,i}.
\end{equation}

We consider the action of the interwiners on
elements of ${\mathcal L}$.  From Proposition $\ref{firstmac}$ it follows that
the coefficient of $\P((-1)^{n-1}a_{-n};q^2, q^4)v_{m,j}$ is
of $P_\lambda(\{q^{-3 +(1/2)}z_j\}_{j=1}^N;q^4, q^{2}).$
We show that the four intertwiners considered so
far leave the lattice ${\mathcal L}$ invariant.

Recall (\cite{Ma}, p.345) that
\begin{equation}\label{skew}
J_{\lambda}(x, z)=\sum_{\mu \subset \lambda}J_{\lambda/\mu}(x)J_{\mu}(z)
\end{equation}
where $x, z$ are infinite sets of indeterminates and the
$J_{\lambda/\mu}(x)$ are the skew integral
Macdonald polynomials.
As before, we will
also consider the restriction to the ring of symmetric functions in a
finite number of indeterminates.

Fix a partition $\lambda$.
Consider the action of $\bpM(z_0)$ on ${\mathcal L}.$  By Proposition
\ref{firstmac}, we have
\begin{equation}
\begin{split} \label{macidentity} &(-q)^{\partial-1} \bpM(z_0)
 (-q)^{\partial-1} \bpM(z_1)
  (-q)^{\partial-1} \bpM(z_2) ...
 (-q)^{\partial-1} \bpM(z_N) \vi  = \\
& \ov{c(z)}\sum_{l(\mu)\leq N}
  J_{\mu}(\{q^{-3+(1/2)}z_j\}_{j=0}^N; q^4, q^2)\mathcal{J}^*_{\mu'}((-1)^{n-1}
a_{-n};q^2, q^4)v_{m, j}
\end{split}
\end{equation}
Now using $\eqref{skew}, \eqref{macidentity}$
\begin{equation}
= \sum_{\ell(\mu) \le N+1} \sum_{\lambda \subset \mu}
J_{\lambda}(\{z_j\}_{j=1}^{N};q^4, q^2)
J_{\mu/\lambda}(z_0;q^4,q^2) {\mathcal J}_{\mu'}^*((-1)^{n-1}a_{-n};
q^2, q^4) v_{m,j}. \notag
\end{equation}
Specializing for a specific $\lambda$ where
$\ell(\lambda) \le N$ we have:
\begin{lem}
\begin{equation}
\begin{align} \notag
& \bpM(z_0) \J^*((-1)^{n-1}a_{-n};q^2, q^4)v_{m,j}   \\
& = \sum_{\mu'; \lambda \subset \mu, \ell(\mu) \le N+1}
\bigl(
 J_{\mu/\lambda}(z_0;q^4,q^2)
\J^*_{\mu'}((-1)^{n-1} a_{-n}; q^2, q^4) v_{m,i}
\bigr). \notag
\end{align}
\end{equation}
\end{lem}

Now from \cite{Ma}, page 340, we have
\begin{equation}
 J_{\mu/\lambda}(z_0;q^4,q^2) = \sum_{\lambda'}
 f_{\mu',\nu'}^{\lambda'} J_{\lambda'},
\end{equation}
where $f_{\mu',\nu'}^{\lambda'}(q^4,q^2) \in
\Z[q,q^{-1}].$
Together with similar calculation for the other cases, this proves
\begin{lem} \label{invarianceone} 
${\mathcal L}$ is invariant under the action of the vertex operators.
\end{lem}

The following lemma is inspired by  \cite{Ma} (see also \cite{Ga}).
\begin{lem}   \label{biglemma}
\begin{equation*}
\prod_{1 \le i < j \le N} (z_i - q^{-2} z_j) =
\sum_{w \in S_N} (-q^{-2})^{\ell(w)} z^{w(\delta)} +
 \sum a_{\gamma_1, \dots, \gamma_n}
 z_1^{\gamma_1}  z_2^{\gamma_2} \dots z_n^{\gamma_n},
\end{equation*}
where $\delta = (N-1, N-2,\dots, 0)$
and for each monomial on the right hand side, some $\gamma_i =
\gamma_j$ for $i \neq j$ and $a_{\ov{\gamma}} \in \Z[q^{-2}],  \ 
a_{\ov{\gamma}}(1) = 0$.
\end{lem}
\begin{pf}
We have
\begin{equation*}
\prod_{i < j} (z_i - q^{-2} z_j) = \sum_{\gamma} (-q^{-2})^{d({\gamma})}
z_1^{\gamma_1} z_2^{\gamma_2} \dots z_N^{\gamma_N},
\end{equation*}
where the summation runs through all $N \times N$ matrices
$(\gamma_{ij})$ of $0$'s and $1$'s such that
\begin{align*} & \gamma_{ii} = 0,\  \gamma_{ij} + \gamma_{ji} = 1 \text{ if
    } i \neq j, \\
      & \text{ and } d(\gamma) = \sum_{i < j} \gamma_{ji}, \ \gamma_i =
      \sum_j \gamma_{ij}.
\end{align*}
When the $\gamma_i$ are all distinct we have
$z_1^{\gamma_1}z_2^{\gamma_2} \dots z_N^{\gamma_N} = z_1^{\gamma_{w(1)}}
\dots z_N^{\gamma_{w(N)}}$ for some permutation $w \in S_N$ and
\begin{equation*}
    \gamma_{w(i)} = \mu_i + (n-i), \ \ (1 \le i \le N),
\end{equation*}
for some partition $\mu: \mu_1 \ge \mu_2 \ge \dots \ge \mu_N \ge 0.$
We claim that all $\mu_i$ are actually $0$, from which it will follow
$z_1^{\gamma_1}
\dots z_N^{\gamma_N} = z^{w(\delta)}$.  In fact, if
$s_{ij} = \gamma_{w(i),w(j)}$,  for $1 \le k \le N$,
\begin{align*}
&  0 \le \mu_1 + \dots + \mu_k  = \sum_{i = 1}^k \sum_{j = 1}^N s_{ij} -
\sum_{i = 1}^k (N-i) \\
& = \frac{1}2 k (k-1) + \sum_{i = 1}^k \sum_{j = k+1}^N s_{ij} - \sum_{i=1}^k
(N-i) \\ & \le \frac{1}2 k (k-1) + k (N-k) - \sum_{i=1}^k(N-i) = 0,
\end{align*}
from which it follows that each $\mu_i = 0.$  Notice that the last
inequality is equality if and only if $s_{ij} = \gamma_{w(i),w(j)} = 1$ for
all pairs $i < j.$  Then, for each distinct $(\gamma_1, \dots \gamma_N)$,
we have:
\begin{equation*} d(S) = \sum_{i<j} s_{ji} = \sum_{i<j}
  \gamma_{w(j),w(i)} = \ell(w),
\end{equation*}
and $d(S)$ is
the number of pairs $i < j$ in $\{1, \dots, N\}$ such that $w(i) >
w(j).$
\end{pf}

\begin{prop}\label{invarianceforphi}
The lattice $\mathcal L$ is invariant under ${x_{k}^\pm}^{(N)}$.  In particular,
$\mathcal L$ is invariant under the action of $U_{\A}.$
\end{prop}

\begin{pf}
By a modification of \eqref{generatingfunction}
the dual integral Macdonald polynomials are also generated by the vertex
operators $\Phi_{\pm}(z)$.  The generators $x_0^\pm, x_{-1}^+, x_{1}^-$ 
satisfy the following commutation relations:
\begin{align*} & [\pP(z), \xom]_{q^{-1}} = 0, && [\pM(z),\xzp] = 0, \\
              & [\pM(z), \xmp] = 0, && [\pP(z), \xom]_q = 0, \\
              & [\pM(z), \xzm] = \pP(z), && [\pP(z),\xzp] = K \pM(z), \\
              & [\pP(z),\xmp] = q^{-1}z^{-1} K^{-1} \pM(z), &&
              [\pM(z),\xom]_{q^{-1}} = q^2 z \pP(z). \\
\end{align*}

By Lemma \ref{invarianceone}, it follows
 that ${\mathcal L}$ is invariant under the $\Phi_\pm(z)$ action, so it
will be sufficient to show that
${x_{k}^\pm}^{(N)} e^{m \alpha} \vi \in {\mathcal L}.$ We compute that
\def\norm{\frac{1}{(2\pi i)^N}}
\begin{align*}
    {x_k^+}^N e^{m \alpha} \vi &= \frac{1}{(2\pi i)^N} \int X^+(z_1)
X^+(z_2) \dots X^+(z_N) z_1^k \dots z_N^k \ dz \  e^{m \alpha} \vi \\
& =  \frac{1}{(2\pi i)^N} \int \exp\bigl(\sum_{n=1}^\infty \frac{(q^n +
  q^{-n})q^{-n/2}}{n} a_{-n} (z_1^n + \dots + z_N^n)\bigr) \\ & \hskip
 .2 in \times \prod_{i < j}
(z_i - z_j)(z_i - q^{-2}z_j) z^{2m + k + i} e^{(m + N)\alpha} \vi dz,
\end{align*}
where we abbreviate $z = z_1 \dots z_N, dz = dz_1 \dots dz_N,$ and the
integration is over the boundary of 
suitable multidimensional disk about $0$.  Observe
that the integrand divided by $\prod_{i<j} (z_i - q^{-2} z_j)$ is an
anti--symmetric function in $z_1, \dots, z_N.$ Invoking Lemma
\ref{biglemma} for $\prod_{i<j} (z_i - q^{-2}z_j)$, we see that
considering antisymmetry, the terms $z_1^{\gamma_1} z_2^{\gamma_2} \dots
z_N^{\gamma_N}$ (for which some $\gamma_i = \gamma_j$) make no contribution
to the integral.  Then
\begin{align*} {x_k^+} & e^{m \alpha} \vi  = \sum_{w \in S_N} 
\int \frac{dz}{(2\pi i)^N} \exp \bigl( \sum_{n=1}^\infty
  \frac{(q^n + q^{-n}) q^{-n/2}}{n}
  a_{-n} (z_1^n + \dots + z_N^n)\bigr) \\ & \times
  \prod_{i<j} (z_i - z_j) (-q)^{-\ell(w)} z^{w(\delta) + (2m + k + i)
    {{\bf 1}}} e^{(m+N) \alpha} \vi \\ & = \bigl( \sum_{w \in S_N}
  q^{-2 \ell(w)} \bigr)
 \int \frac{dz}{(2\pi i)^N} \exp\bigl( \sum_{n = 1}^\infty
  \frac{(q^n + q^{-n}) q^{-n/2}}{n} a_{-n} (z_1^n + \dots + z_N^n)
  \bigr) \\ & \times \prod_{i < j} (z_i - z_j) z^{\delta + (2m + k +
    i){{\bf 1}}} e^{(m + N) \alpha} \vi,
\end{align*}
where $\delta = (N-1, \dots, 0)$ and ${{\bf 1}} = (1, 1, \dots, 1).$

By the orthogonality of $J_\lambda$ \eqref{innerproductPi} we have
\begin{equation*}
\begin{split}
\exp\bigl( \sum_{n=1}^\infty \frac{(q^n + q^{-n})q^{-n/2}}{n}
      a_{-n} & (z_1^n + \dots + z_N^n)\bigr) = \\ &
     \sum_{\ell(\lambda) \le N} {\mathcal J}^*_\lambda(a_{-n}; q^2, q^4)
     J_\lambda(z_i;q^2, q^4) q^{-|\lambda|}.
\end{split}
\end{equation*}
The integrality of $J_\lambda$ (see for example \cite{GT}) implies that
\begin{equation*}
  J_{\lambda}(z_i; q^2, q^4) = \sum_{\lambda \ge \mu} a_{\lambda \mu}
  s_{\mu}(z_i) , \ a_{\lambda \mu} \in \Z[q].
\end{equation*}
Then, since $\sum_{w \in S_N} q^{-2\ell(w)} = q^{-\binom{N}{2}} [N]!,$
we have
\begin{align*} {x_k^+}^{(N)} e^{m \alpha} \vi & = q^{- \binom{N}{2}}
  \sum_{\ell(\lambda) \le N} {\mathcal J}_\lambda^*(a_{-n};q^2, q^4) \int
\frac{dz}{(2\pi i)^N}
  J_\lambda(z_i;q^2,q^4) q^{-|\lambda|} \\ & \times \prod_{i<j} (z_i - z_j)
  z^{\delta + (2m + k +i) {\bf 1}} e^{(m+N)\alpha} \vi
   \\
   & = q^{-\binom{N}{2}} \sum_{\ell(\lambda) \le N}
   {\mathcal J}^*_\lambda(a_{-n};q^2,q^4) q^{-|\lambda|} e^{(m+N) \alpha} \vi \cdot
   \sum_{\lambda \ge \mu} a_{\lambda \mu} 
\\ & \times \int \frac{dz}{(2\pi i)^N}  
        s_\mu(z_i) 
    \prod_{i < j}
   (z_i - z_j) z^{\delta + (2m + k + i) {\bf 1}},
 \end{align*}
where the last integral is an integer by the integrality of Schur
functions.  Therefore it follows ${x_k^+}^{(N)}
e^{m \alpha} \vi \in {\mathcal L}.$ The case of ${x_k^-}^{(N)}$ can be proved
similarly.
\end{pf}


\section{Canonical basis and Macdonald Polynomials}

\begin{defn} Let $m \in \Z,\ i = 0,1.$ Define $A_{(n,m,i)}$ to 
be the subspace of $\H$ spanned by 
$\tilde {\mathcal Q}_\lambda v_{m,i},$ where $|\lambda| = n.$
\end{defn}
Clearly $\H = \oplus_{n,m,i}  \Q(q) A_{(n, m, i)}$ and
the canonical basis respects this grading.
\begin{prop} \label{sector} Let $b \in \H_\A$ be an element of the
  canonical basis of $\H$ in $A{(n,m,i)}.$ Write $b = p(a_{-k}) \otimes
  v_{m,i}$ where $p$ is a polynomial of degree $n$.  Then for any $m'$,
  up to a sign, $b' = p(a_{-k})\otimes v_{m',i}$ is also an element of
  the canonical basis of $\H.$
\end{prop}
\begin{pf} This follows from the coincidence of Macdonald's and
  Kashiwara's forms.  By the characterization of the canonical basis,
  $b = p(a_{-k}) \otimes v_{m,i}$ is in $U_{\mathcal A}$, bar--invariant,
  and $(b,b) = 1 + q^{-1} f(q^{-1}).$  The same holds for
$b' = p(a_{-k})\otimes v_{m',i}.$
\end{pf}

 The
integrality result (see, for example, \cite{GT}) for the two variable
Kostka matrix $K(q,t)$ implies:
\begin{equation}
J_\lambda(q,t) = \sum_{\mu \le \lambda, |\mu| = |\lambda|}
v_{\lambda,\mu}(q,t) m_{\mu}, \ \ v_{\lambda,\mu} \in \Z[q],
\end{equation}
where the  $m_\mu$ are the
monomial symmetric functions.

Combining this with
\begin{equation} s_\lambda = \sum_{\mu \le \lambda} K_{\lambda,\mu}
  m_{\mu},
\end{equation}
where $(K_{\lambda, \mu})$
is the usual Kostka matrix, and $s_\lambda$ are the Schur functions,
we see that
\begin{equation} \label{integralM}
 J_{\lambda}(q,t) = \sum_{\mu \le \lambda} w_{\lambda,\mu}(q,t)
s_\mu, \text{ where } w_{\lambda,\mu} \in \Z[q,t].
\end{equation}

\begin{defn} Let $m \in \Z,\ i = 0,1.$
Let the Schur polynomials
$\tilde s_\lambda v_{m,i} = \tilde \P_\lambda
(\hat a_{-n};$ $q ,q)v_{m,i}.$
\end{defn}

By results of \cite{CP} it is known that the $\tilde s_{\lambda}
v_{m,i}$ are contained in the lattice of divided powers, and it follows
from \eqref{integralM} that $\tilde \J_\lambda$ is contained in the
lattice of divided powers.

There is a natural order (see \cite{LTT}) on the canonical basis.
Since the canonical basis respects the grading $A_{(n, m, i)}$ of $\H$,
we can consider in each graded component the transition matrix between
the dual canonical basis and the integral Macdonald polynomials, given by:
\begin{equation}    \tilde {\mathcal J}_\lambda v_{m,i}
 = \sum_{\mu, \ \sum \mu_i = n}
  a_{\mu,\lambda}(q) B^*_\mu
\end{equation}
where $a_{\mu, \lambda} \in \Z[q,q^{-1}]$,
and $B^*_\mu$ are elements of the dual canonical basis in $A_{n,m,i}$.

Let $C$ be the diagonal matrix consisting of $c_\lambda(q^4,q^2)$.
\begin{prop} Let $A = C^{-1}(a_{\mu,\lambda}).$
\begin{enumerate}
\item The matrix $A$ consists of bar invariant elements.
\item  A consists of polynomials in $q$ and $q^{-1}$ with integral
  coefficients.
\end{enumerate}
\end{prop}
\begin{pf}
1) follows from the previous section. 2) follows from the discussion
above.
\end{pf}

From Proposition \ref{orthogonal}, we see that the polynomials $\tilde
\J_\lambda v_{m,i}$ form a quasi--orthogonal basis of $H$.  However,
they are not elements of the canonical basis except in the case where
$\lambda$ is the empty partition.

\begin{conj} The matrix $A$ is upper unitriangular with coefficients
in $\N[q,q^{-1}].$
\end{conj}

Fix a weight and sector as in Proposition \ref{sector}. With respect to
this $A_{(n,m,i)}$ let
$B$ (resp. $B^*$) denote the canonical basis (dual canonical basis) of
$\H$. Let $J$ (resp. $J^*$) denote the basis $\tilde \J_\lambda$
($\tilde \J^*_\lambda$) for this sector.  Let 
$L$ denotes the diagonal matrix consisting of $c_\lambda(q^4,q^2)
c'_\lambda(q^4,q^2).$  Then:
$$ B = A(q)^{t} J^* =  A(q)^t L^{-1} J = A(q)^t
L^{-1} A(q) B^*.$$

We rewrite this in terms of the $\tilde \P$ and $\tilde {\mathcal Q}$.
Since $b_\lambda^{-1}(q^2,q^4)\tilde {\mathcal Q}_\lambda(q^2,q^4) =
\tilde \P_\lambda(q^2,q^4)$, we have immediately that the basis of Macdonald
polynomials is a ``square root'' of the transition matrix from the dual
canonical basis to the canonical basis.
$$  B = A(q)^{t} D(q) A(q) B^*,$$
where $D(q)$  is the  diagonal matrix diag$(b^{-1}_\lambda(q^2,q^4))$.

We conclude the paper with a brief discussion of the
possible nature of the coefficients of the transition matrix $A(q)$ that
also suggests an approach and provides support for our conjecture.  It
was suggested in \cite{DFJMN} that one can define an imbedding
of $\H = V(\Lzero) \oplus V(\Lone)$ into the infinite product $V_1
\otimes V_1 \otimes ...$ of two dimensional representations of
$U_q({\frak sl}_2)$ via the correspondence of the canonical bases of
both spaces.  On the other hand, the dual canonical basis in a finite
product ${V_1}^{\otimes n}$
of the two dimensional representations of $U_q({\frak sl}_2)$ has been
studied recently in \cite{FK2}.   In the finite product case, the
transition matrix  from the dual canonical basis to the canonical basis also
admits a factorization
\begin{equation} B_n = A_n(q)^t A_n(q) B_n^* \notag
\end{equation}
via the elementary basis consisting of the elements $\{v_{\varepsilon_1}
\otimes \dots \otimes v_{\varepsilon_n}, \varepsilon_i = \pm \}.$ Direct
combinatorial arguments show that the matrix $A_n(q)$ is upper
unitriangular with coefficients in $\N[q,q^{-1}].$ A further study of
the transition matrix $A_n(q)$ in \cite{FKK} shows that its coefficients
are identified with the Kazhdan--Lusztig polynomials for Grassmanians,
i.e. those corresponding to maximal parabolic subgroups in $S_n.$ It
turns out that precisely in this case there exist simple combinatorial
formulas for Kazhdan--Lusztig polynomials first obtained in \cite{LS}
and rederived via the graphical calculus for $U_q({\frak sl}_2)$ in
\cite{FKK}.  Moreover, one can invoke the representation theoretic
interpretation of the Kazhdan--Lusztig polynomials associated to $S_n$
as the Jordan--Holder multiplicities of irreducible representations in
Verma modules for ${\frak sl}_n$, and this yields the  desired properties of
$A_n(q)$ without explicit calculation.  To obtain a similar
interpretation of the transition matrix from the canonical basis in $\H$
to Macdonald polynomials one has to make sense of the limit of the above
finite dimensional construction when $n \rightarrow \infty,$ which is a
very delicate matter.  However, when it is done correctly, the explicit
formulas for the coefficients of $A(q)$ should be even more elementary
than the ones for $A_n(q).$ This provides some assurance that the symmetric
functions corresponding to the canonical basis have a simple enough
description relative to Macdonald polynomials.  At the present moment,
however, we do not know to which symmetric functions they correspond.


\bigskip

\noindent{\bf Acknowledgements.} We would like to thank H. Garland and
I. Grojnowski for their active interest and participation at
different stages of this work. I. Frenkel is supported in part by NSF grants
DMS--9400908, DMS--9700765.  
N. Jing is supported in part by NSA grant MDA 904-97-1-0062.

\end{document}